\documentclass[12pt]{article}

\usepackage[latin1]{inputenc}
\usepackage{bbm}
\usepackage{stmaryrd}
\usepackage{latexsym}
\usepackage{amsfonts}
\usepackage{amsmath,amssymb,amscd}
\usepackage{yfonts}
\usepackage{dsfont}
\usepackage{mathabx}
\usepackage[dvips]{graphicx}
\usepackage{epsfig}
\usepackage{psfrag}

\DeclareFontFamily{U}{txsyc}{}
\DeclareFontShape{U}{txsyc}{m}{n}{
   <-> txsyc%
}{}
\DeclareFontShape{U}{txsyc}{bx}{n}{
   <-> txbsyc%
}{}
\DeclareFontShape{U}{txsyc}{l}{n}{<->ssub * txsyc/m/n}{}
\DeclareFontShape{U}{txsyc}{b}{n}{<->ssub * txsyc/bx/n}{}
\DeclareSymbolFont{symbolsC}{U}{txsyc}{m}{n}
\SetSymbolFont{symbolsC}{bold}{U}{txsyc}{bx}{n}
\DeclareFontSubstitution{U}{txsyc}{m}{n}
\DeclareMathSymbol{\df}{\mathrel}{symbolsC}{"42}
\DeclareMathSymbol{\fd}{\mathrel}{symbolsC}{"43}
\DeclareMathSymbol{\lJoin}{\mathrel}{symbolsC}{"58}
\DeclareMathSymbol{\rJoin}{\mathrel}{symbolsC}{"59}

\newcommand{\cB}{{\cal B}}
\newcommand{\cC}{{\cal C}}

\newcommand{\cF}{{\cal F}}
\newcommand{\cG}{{\cal G}}
\newcommand{\cH}{{\cal H}}

\newcommand{\cL}{{\cal L}}

\newcommand{\cR}{{\cal R}}

\newcommand{\DD}{\mathbb{D}}
\newcommand{\EE}{\mathbb{E}}

\newcommand{\LL}{\mathbb{L}}
\newcommand{\NN}{\mathbb{N}}
\newcommand{\PP}{\mathbb{P}}

\newcommand{\RR}{\mathbb{R}}

\newcommand{\ZZ}{\mathbb{Z}}

\newcommand{\cCb}{{\cal C}_{\mathrm{b}}}
\newcommand{\cCc}{{\cal C}_{\mathrm{c}}}

\newcommand{\iy}{\infty}
\newcommand{\lt}{\left}
\newcommand{\me}{\medskip}

\newcommand{\pa}{\partial}
\newcommand{\ri}{\rightarrow}
\newcommand{\rt}{\right}
\newcommand{\sm}{\smallskip}
\newcommand{\tr}{\triangle}

\newcommand{\wi}{\widetilde}
\newcommand{\wit}{\widehat}

\newcommand{\fo}{\forall\ }

\newcommand{\lve}{\lt\vert}
\newcommand{\lVe}{\lt\Vert}

\newcommand{\rve}{\rt\vert}
\newcommand{\rVe}{\rt\Vert}

\newcommand{\st}{\,:\,}

\newcommand{\un}{\mathds{1}}

\newcommand{\bq}{\begin{eqnarray*}}
\newcommand{\bqn}[1]{\begin{eqnarray}\label{#1}}
\newcommand{\eq}{\end{eqnarray*}}
\newcommand{\eqn}{\end{eqnarray}}
\newcommand{\wwtbp}{\par\hfill $\blacksquare$\par\me\noindent}
\newcommand{\thistitlepagestyle}{}

\newcommand{\ttsim}{\raise.17ex\hbox{$\scriptstyle\mathtt{\sim}$}}

\newtheorem{pro}{Proposition} 
\newtheorem{cor}[pro]{Corollary}
\newtheorem{lem}[pro]{Lemma}
\newtheorem{theo}[pro]{Theorem}
\renewcommand{\thepro}{\arabic{pro}}

\newenvironment{rem}
{\par\me\refstepcounter{pro}\noindent{\bf Remark \thepro\ }}
{\par\hfill $\square$\par\me\noindent}
\newcommand{\proof}{\par\me\noindent\textbf{Proof}\par\sm\noindent}

\title{Strong stationary times for one-dimensional diffusions}
\author{Laurent Miclo
}
\date{\box1
}
\newcommand{\rE}{\mathring{E}^*}

\begin{document}

\setbox1=\vbox{
\large
\begin{center}
Institut de Mathématiques de Toulouse, UMR 5219\\
Université de Toulouse and CNRS, France\\
\end{center}
} 
\setbox3=\vbox{
\hbox{miclo@math.univ-toulouse.fr\\}
\vskip1mm
\hbox{Institut de Mathématiques de Toulouse\\}
\hbox{Université Paul Sabatier\\}
\hbox{118, route de Narbonne\\} 
\hbox{31062 Toulouse Cedex 9, France\\}
}
\setbox5=\vbox{
\box3
}

\maketitle
\thistitlepagestyle
\abstract{
A necessary and sufficient condition is obtained for the existence of strong stationary times
for ergodic one-dimensional diffusions, whatever the initial distribution.
The strong stationary times are constructed through intertwinings with dual processes, in the Diaconis-Fill sense, taking values in the set
of segments of the extended line $\mathbb{R}\sqcup\{-\infty,+\infty\}$. They  can be seen as natural $h$-transforms of the extensions to the diffusion framework of the
evolving sets of Morris-Peres. Starting from a singleton set, the dual process begins by evolving into true segments
in the same way a Bessel process of dimension 3 escapes from 0. The strong stationary time corresponds to the first
time the full segment $[-\infty,+\infty]$ is reached.
The benchmark Ornstein-Uhlenbeck process cannot be treated in this way, it will nevertheless be seen how to use 
other strong times to recover its optimal exponential rate of convergence in the total variation sense.
}
\vfill\null
{\small
}\par

\newpage 

\section{Introduction}

A strong stationary time $\tau$ is a stopping time relative to the filtration generated by an ergodic 
Markov process $(X_t)_{t\geq 0}$ (and possibly some independent randomness) which is such that $\tau$ and $X_\tau$ are independent
and $X_\tau$ is distributed according to the underlying invariant probability distribution.
They were first introduced by Aldous and Diaconis \cite{MR876954} in the context of finite Markov chains.
Staying in the finite framework, Diaconis and Fill \cite{MR1071805} developed the important tool
of intertwining with absorbed Markov chains to construct strong stationary times.
Intertwining of diffusions was also investigated by 
Rogers and Pitman \cite{MR624684} and
Carmona, Petit and Yor  \cite{MR1654531},
especially 
to deduce identities in law for particular  processes.
Recently, Pal and Shkolnikov \cite{2013arXiv1306.0857P} studied some conditions insuring that there
exists an intertwining between two Markov semi-groups and their
article also provides a welcome survey of applications of intertwining relations.
Our goal here is to come back to the investigation of strong stationary times 
through intertwining, but in the context of diffusions. We will also point out a  
relation with an extension to this continuous setting of the evolving sets of Morris and Peres \cite{MR2198701}.
More precisely, we are to be mainly concerned with one-dimensional diffusions,
the simplest continuous framework and nevertheless already displaying some interesting features.
Of course, extensions to multidimensional situations are more promising and challenging,
there are outside the scope of this paper, which can be seen as only working out the preliminary steps 
in this direction that we hope to 
investigate 
in the future.\par\me
Consider the one-dimensional Markov generator given by
\bqn{L}
L&\df& a\pa^2+b\pa\eqn
where $a>0$ and $b$ are two functions defined on $\RR$. 
We won't be interested in regularity issues, so we assume
that they are smooth and $L$ can be interpreted as an operator from $\cCc^{\infty}(\RR)$ to itself.
It is often convenient to extend it as a self-adjoint operator on a $\LL^2$ space.
Indeed, consider
\bqn{mu}
\nonumber\fo x\in\RR,\qquad c(x)&\df& \int_0^x\frac{b(y)}{a(y)}\, dy\\
\mu(x)&\df& \frac{\exp(c(x))}{a(x)}
\eqn
We will denote by the same symbol $\mu$ the measure admitting the function $\mu$ as density 
with respect to the Lebesgue measure.
It is well-known (cf.\ for instance 
the chapter 15 of the book of Karlin and Taylor \cite{MR611513}), and elementary to recover, that
the operator $L$ is symmetric in $\LL^2(\mu)$, so we 
can consider the corresponding Freidrich's extension.
\par\sm
Since we are 
only interested in positive recurrent diffusions, we begin by making the assumption that $\mu$ is a finite measure:
\bqn{mfini}
m\ \df\ \int_\RR  \frac{\exp(c(x))}{a(x)}
\, dx &<&+\infty\eqn
and we renormalize $\mu$ into a probability distribution, 
replacing (\ref{mu}) by
\bqn{mu2}
\mu(x)&\df& \frac{\exp(c(x))}{m\,a(x)}\eqn
Let $X\df(X_t)_{t\geq 0}$ be a diffusion process whose generator is $L$.
The above finiteness assumption
does't prevent $X$ from exploding in finite time.
Indeed the general criterion for $X$ to be non-explosive is that 
 \bq
\int_{-\iy}^0\mu([y,0])\exp(-c(y))\,dy\ =\ +\iy&\hbox{ and }&\int_0^{\iy}\mu([0,y])\exp(-c(y))\,dy\ =\ +\iy\eq
(see for instance Theorem 3.2 (3) of Chapter 6 of the book of Ikeda and Watanabe \cite{MR1011252}).
When $\mu$ is finite, as it is implicit throughout the paper except otherwise stated, this condition reduces to
 \bqn{JJ}
\int_{-\iy}^0\exp(-c(y))\,dy\ =\ +\iy&\hbox{ and }&
\int_0^{\iy}\exp(-c(y))\,dy\ =\ +\iy\eqn
A diffusion $X$ whose generator $L$ satisfies \eqref{mfini} and \eqref{JJ} is said to be positive recurrent.
\par\sm
The process $X$ is a priori defined on a probability space $(\Omega,\cF,\PP)$ endowed with the filtration $(\cF_t)_{t\geq 0}$
generated by $X$. For instance, $\Omega$ can be taken to be the set of continuous trajectories $\cC([0,+\iy),[-\iy,+\iy])$
endowed with the $\sigma$-field and the filtration generated by the canonical coordinate process.
But to allow for extra randomness, it is useful to enlarge the initial setting $(\Omega,\cF,\PP,(\cF_t)_{t\geq 0})$ 
into  $(\bar\Omega,\bar\cF,\bar\PP,(\bar \cF_t)_{t\geq 0})$, preserving the fact that $X\df(X_t)_{t\geq 0}$ is a continuous process 
starting from $x_0$, 
Markovian with respect to the filtration $(\bar\cF_t)_{t\geq 0}$
and whose generator is $L$. This is often done
by considering the tensor product of $(\Omega,\cF,\PP,(\cF_t)_{t\geq 0})$ with another probability space.
\par
A random time $\tau$ taking values in $[0,\iy]$ is said to be a 
stopping time,  if it is defined on a framework  $(\bar\Omega,\bar\cF,\bar\PP,(\bar \cF_t)_{t\geq 0})$ as above   and if it 
a stopping time with respect to the filtration $(\bar\cF_t)_{t\geq 0}$, namely if
\bq
\fo t\geq 0,\qquad \{\tau\leq t\} \in\bar\cF_t\eq
From a practical point of view, it means that $\tau$ is constructed from $X$ and from some independent randomness $Y$ in such a way that for any $t\geq 0$ and 
in view of $Y$, to decide whether $\tau\leq t$  or not, it is sufficient to look at the trajectory $X_{[0,t]}\df (X_s)_{s\in[0,t]}$.\par
The stopping time $\tau$, taking values in $[0,\iy)$, is said to be strong, if $\tau$ and $X_{\tau}$ are independent.
It is said to be a strong stationary time, if  furthermore $X_\tau$ is distributed according to $\mu$.
\par\sm
Our main goal in this paper is to investigate the existence of strong stationary times for $X$.
To state our first result, we need the following quantities
\bq
I_-&\df&\int_{-\iy}^0\lt(\int_x^0
\exp(-c(y))\, dy \rt) 
\,\mu(dx)\\
I_+&\df&\int_0^{+\iy}\lt(\int_0^{x}
\exp(-c(y))\, dy\rt) 
\,\mu(dx)\\
I&\df& \max(I_-, I_+)\eq
Only the finiteness of $I$ will be important for us and for that the role of 0 is irrelevant: it could be replaced
by any other point of $\RR$.
But if we were looking for quantitative bounds, it should be chosen more carefully, maybe replacing it by $x_0$
in the case where $X$ starts from the initial deterministic condition $X_0=x_0$.
For the next result, we allow any initial distribution for $\cL(X_0)$.
\par
\begin{theo}\label{t1}
Assume that $X$ is positive recurrent.
There exists a strong stationary time for $X$, whatever its initial distribution, if and only if $I<+\iy$.
\end{theo}
\par
\begin{rem}
Despite we made in this paper the deliberate choice not to get involved in optimal regularity questions, let us
mention that the natural framework for the previous result is that of general one-dimensional diffusions (see for instance 
Section 5.3 of the book of Revuz and Yor \cite{MR1725357}):
the generator is no longer described by (\ref{L}), but under the form 
$L=\frac{d}{dm}\frac{d}{ds}$, where $m$ is the speed measure 
and $s$ is the scale function
(in our regular setting, they admit densities with respect to the Lebesgue measure  respectively given
by $\exp(c)/a$ and $\exp(-c)$).
 In this context and up to a constant factor, 
the quantity analogous to $I$ writes down as $\max(\bar I_-,\bar I_+)$, where
\bq
\bar I_-\ \df\ \int_{-\iy}^0s([x,0])
\,m(dx)
&\hbox{ , }& \bar I_+\ \df\ \int_0^{+\iy}s([0,x])\, m(dx)
\eq
and we expect Theorem \ref{t1} to be still true. Remark that the positive recurrence of $X$ can also be expressed
through $m$ and $s$ only: $m(\RR)<+\iy$ and $s(\RR_-)=s(\RR_+)=+\iy$.
\end{rem}
\begin{rem}\label{Mao}
Instead of the whole line $\RR$, we could have considered the half-line $\RR_+$ with usual reflection at 0.
Similar notions can be introduced in this context and the arguments can be adapted to show that
the corresponding Theorem \ref{t1} is valid, where $I$ is replaced by $I_+$.
In a recent preprint, Cheng and Mao \cite{Cheng_eigentime} showed that the assumption $I_+<+\iy$
is equivalent to several conditions, among which the strong ergodicity of $X$ and the fact that the essential spectrum of $L$
is empty and that the sum of the inverses of its non-zero eigenvalues is finite.
This amounts to say that the associated centered Green operator (which is the inverse of the generator on the space
of functions whose mean with respect to the invariant measure vanishes) has a finite trace.
While this result was only stated for half-space, we strongly believe it also holds for ergodic diffusion on $\RR$.
Thus  the existence of strong stationary times of a recurrent positive one-dimensional diffusion, whatever the initial condition, would be equivalent to its centered Green operator having a finite trace.
At least, this is coherent with the fact for finite birth and death chains starting from a boundary, the optimal strong stationary time
is distributed as a sum of independent exponential variables of parameters the inverses of the absolute values of the non-zero eigenvalues 
of the generator (see Fill \cite{MR1144727}).
This result has been recently extended by Cheng and Mao \cite{Cheng_passage} to diffusions on a  compact segment of $\RR$ with reflecting boundaries, when the process starts from one of those boundaries. 
\par
We believe that the finiteness of the trace of the centered Green operator is always a sufficient condition for the existence of strong stationary times,
but the necessity 
of this property cannot be true in full generality: consider a probability $\mu$ on  a general measurable space and
let $L$ be the generator acting on functions $f\in\LL^2(\mu)$ by $L[f]\df\mu[f]\un-f$.
A strong stationary time, whatever the initial distribution, is given by the first jump. The spectrum of $L$ consists of 0 (with multiplicity 1)
and of 1  with multiplicity the dimension of $\{f\in\LL^2(\mu)\st \mu[f]=0\}$. So if the latter dimension is infinite, we get a counter-example to the necessity condition outside the framework of one-dimensional diffusions.
\end{rem}
As announced at the beginning of this section, a strong stationary time will be constructed through duality via 
intertwining relations. More precisely, 
let
\bq
E^*&\df& \{(x,y)\st x,y\in[-\iy,+\iy],\, x\leq y\}\setminus\{(-\iy,-\iy),(+\iy,+\iy)\}\\
\rE&\df& \{(x,y)\in\RR^2\st x<y\}\eq
be the interior of $E^*$ and $D^*\df\{(x,x)\st x\in\RR\}\subset E^*$ be the diagonal of  $\RR^2$.
Consider the Markov kernel $\Lambda$ from $E^*$ to $\RR$ defined by
\bq
\fo (x,y)\in E^*,\,\fo A\in\cB(\RR),\qquad
\Lambda((x,y),A)&\df&\lt\{\begin{array}{ll}
\delta_{x}(A)&\hbox{, if $y=x$}\\
\noalign{\vskip 2mm}
\frac{\mu([x,y]\cap A)}{\mu([x,y])}
&\hbox{, otherwise}
\end{array}\rt.\eq
where $\cB(\RR)$ stands for the set of Borelian sets from $\RR$.
\par
Transposing to the diffusion setting the program described by Diaconis and Fill \cite{MR876954} for
finite Markov chains, we are looking for a diffusion generator $L^*$ on $E^*$ satisfying
the intertwining relation $\Lambda L=L^* \Lambda$, in the sense
that at least on $E^*\setminus( D^*\sqcup \{(-\iy,+\iy)\})$,
\bqn{inter}
\fo f\in\cCc^{\iy}(\RR),\qquad \Lambda[ L[f]]&=&L^* [\Lambda[f]]\eqn
Here is one solution: on $\rE$,
\bqn{Lstar}
 L^*&\df&
(\sqrt{a(y)}\pa_y-\sqrt{a(x)}\pa_x)^2+(a'(x)/2-b(x))\pa_x+(a'(y)/2-b(y))\pa_y\\&&\nonumber+2\frac{\sqrt{a(x)}\mu(x)+\sqrt{a(y)}\mu(y)}{\mu([x,y])}(\sqrt{a(y)}\pa_y-\sqrt{a(x)}\pa_x)
\eqn
while on $\RR\times\{+\iy\}$,
\bqn{Lstar2}
 L^*&\df&
(\sqrt{a(x)}\pa_x)^2+(a'(x)/2-b(x))\pa_x-2\frac{\sqrt{a(x)}\mu(x)}{\mu([x,+\iy))}\sqrt{a(x)}\pa_x
\eqn
and on $\{-\iy\}\times\RR$,
\bqn{Lstar3}
L^*&\df&
(\sqrt{a(y)}\pa_y)^2+(a'(y)/2-b(y))\pa_y+2\frac{\sqrt{a(y)}\mu(y)}{\mu((-\iy,y])}\sqrt{a(y)}\pa_y
\eqn
Formally, (\ref{Lstar2}) and (\ref{Lstar3}) are obtained by respectively replacing $y$ by $+\iy$ and $x$ by $-\iy$
in (\ref{Lstar}). Such extensions of (\ref{Lstar}) will be called natural in the sequel.\par
We put a Dirichlet condition at $(-\iy,+\iy)$, insuring that it is an absorbing point.\par
It is not necessary to make precise the boundary condition on the diagonal $D^*$, 
because it is an entrance boundary:
\par
\begin{pro}\label{entrance}
For any $x_0\in\RR$, there is a continuous Markov process $Z^*\df(Z^*_t)_{t\geq 0}$ starting from $(x_0,x_0)$, whose generator is $L^*$ (in the sense of martingale problems) and satisfying for all $t>0$, $Z^*_t\in E^*\setminus D^*$.
The law of this process is unique if we impose that after the possibly finite time
\bqn{taustar}
\tau^*&\df& \inf\{t\geq 0\st Z^*_t=(-\iy,+\iy)\}\eqn
$Z^*$ stays at  position $(-\iy,+\iy)$ (i.e.\ if we consider
the minimal process). 
\end{pro}\par
The generator $L^*$ defined in (\ref{Lstar}) is not the unique one satisfying (\ref{inter}).
This relation is also true if $L^*$ is replaced by
\bqn{dagger}
\check L^*&\df&(\sqrt{a(y)}\pa_y+\sqrt{a(x)}\pa_x)^2+(a'(x)/2-b(x))\pa_x+(a'(y)/2-b(y))\pa_y\\&&\nonumber+2\frac{\sqrt{a(y)}\mu(y)-\sqrt{a(x)}\mu(x)}{\mu([x,y])}(\sqrt{a(y)}\pa_y+\sqrt{a(x)}\pa_x)\eqn
(on $\rE$ and its natural extensions on $\RR\times\{+\iy\}$ and $\{-\iy\}\times\RR$).
For this operator, $D^*$ is not an entrance boundary: an associated process starting on $D^*$ stays in $D^*$, this is related  to the fact
that the mapping $E^*\setminus D^*\ni(x,y)\mapsto (\sqrt{a(y)}\mu(y)-\sqrt{a(x)}\mu(x))/\mu([x,y])$ can be naturally  extended into
a symmetric and smooth function on $\RR^2$.
\\
There are other generators satisfying (\ref{inter}), e.g.\ 
 the elliptic operator
\bqn{ddagger}
 \breve L^*&\df&a(y)\pa^2_y+a(x)\pa^2_y+(a'(x)/2-b(x))\pa_x+(a'(y)/2-b(y))\pa_y\\&&\nonumber+2\frac{1}{\mu([x,y])}(a(y)\mu(y)\pa_y-a(x)\mu(x)\pa_x)\eqn
 (on $\rE$ and its natural extensions on $\RR\times\{+\iy\}$ and $\{-\iy\}\times\RR$).
One would have remarked that $\breve L^*=(L^*+\check L^*)/2$ and more generally for any $\alpha\in(0,1)$, the generator
$(1-\alpha) L^*+\alpha \check L^*$ satisfies (\ref{inter}) and is elliptic. But as it will be seen in Remark \ref{alpha}
at the end of the next section, these generators lead to strong stationary times which are larger than those obtained from $L^*$.
\par\sm
The generator $L^*$ defined in (\ref{Lstar}) has another interest: it is related via a Doob transform to the continuous equivalent of 
the evolving sets introduced by Morris and Peres \cite{MR2198701} for denumerable Markov chains.
Consider the generator given by
\bqn{wiL}
\wi L&\df&
(\sqrt{a(y)}\pa_y-\sqrt{a(x)}\pa_x)^2+(a'(x)/2-b(x))\pa_x+(a'(y)/2-b(y))\pa_y\eqn
(on $\rE$ and its natural extensions on $\RR\times\{+\iy\}$ and $\{-\iy\}\times\RR$).
The diagonal is not an entrance boundary for this generator and to associate a (minimal) Markov process $\wi Z\df(\wi X_t,\wi Y_t)_{t\geq 0}$,
we can impose Neumann boundary condition on $D^*$. It amounts to see $\wi L$ as a generator on $\RR^2$ and for any point $(x,y)$ of the plane,
to identify  $(x,y)$ with $(y,x)$.
The segment valued process $([\wi X_t,\wi Y_t])_{t\geq 0}$ is then a continuous evolving set in $\RR$.
Next consider the mapping $h$ defined on $E^*$
by
\bqn{h}
\fo z=(x,y)\in E^*,\qquad h(z)&\df&\mu([x,y])\eqn
It will be checked in Lemma \ref{evolset} of next section that $\wi L[h]=0$ on $\rE$.
Then $L^*$ is the Doob transform of $\wi L$ through $h$:
\bq
L^*[\cdot]&=&\frac1h\wi L[h\,\cdot\,]\\
&=&\wi L[\cdot]+\wi \Gamma[\ln(h),\cdot]\eq
where $\wi \Gamma$ is the carré du champ associated to $\wi L$:
for any smooth functions $f,g$ defined on $\rE$, 
\bq
\wi \Gamma[f,g]&\df& \wi L[fg]-f \wi L[g]-g \wi L[f]\eq
\par\me
Let us now come back to a diffusion process $X$ 
as in Theorem \ref{t1} and denote by $m_0$ its initial distribution.
Consider the probability $m^*_0$ defined on $E^*$ by $m_0^*\df \int \delta_{(x,x)}\, m_0(dx)$,
so that $m_0^*\Lambda =m_0$.
In general it is not the only probability on $E^*$ satisfying this relation, for instance if $m_0=\Lambda(z,\cdot)$,
with $z\in \rE$, it seems more appropriate to choose $m_0^*\df \delta_z$.
The strong stationary time constructed in Proposition \ref{Itau} below does depend on the choice of $m_0^*$,
but in this paper we will not consider the important question of finding the best possible choice for $m_0^*$
(next section will show how to construct a process $Z^*$ starting from any initial distribution on $E^*$, indeed Proposition \ref{entrance} presented the most difficult cases).
As it is explained by 
Diaconis and Fill \cite{MR1071805} in the finite setting,
the relations $m_0^*\Lambda =m_0$ and (\ref{inter}) should enable to couple
$X$ with the process $Z^*$, defined similarly as in Proposition \ref{entrance} but with $m_0^*$ as initial distribution, in such a way that for any $t\geq 0$,
the conditional law of $X_t$ knowing the trajectory $Z^*_{[0,t]}$ is given by $\Lambda(Z^*_{t},\cdot)$.
The extension to  the positive recurrent one-dimensional diffusion case turned out to be quite tricky  and will be developed in Section \ref{SecI}
% for our setting
(unfortunately the results of
Pal and Shkolnikov \cite{2013arXiv1306.0857P} cannot be applied 
straightforwardly).
Let us admit this technical point for the time being. A convenient feature of this coupling is that
it can be obtained by starting with a trajectory $X$ and by constructing 
$Z^*$ from $X$ and independent randomness. More precisely, for any $t\geq 0$, the  piece of trajectory $Z^*_{[0,t]}$
is constructed from $X_{[0,t]}$ and independent randomness. Thus any stopping time $\tau$
with respect to the filtration generated by
the process $Z^*$ is also a stopping time for $X$.
This is important, because the
previous conditional property extends to any finite stopping time $\tau$ with respect to the filtration generated by
the process $Z^*$:
\bqn{condtaustar}
\cL(X_{\tau}\vert Z^*_{[0,\tau]})&=& \Lambda(Z^*_{\tau},\cdot)\eqn
where the l.h.s.\ is the conditional law of $X_{\tau}$ knowing the trajectory
$Z^*_{[0,\tau]}$.\\
In particular if we consider the stopping time $\tau^*$ defined in (\ref{taustar})
and if we impose conditions such that this $Z^*$-stopping time is a.s.\ finite, then
it is a strong stationary time for $X$.
Indeed, the above considerations show that $\tau^*$ is a stopping time for $X$.
Next note that $X_{\tau^*}$ is independent from $Z^*_{[0,\tau^*]}$, because
according to (\ref{condtaustar}), $\cL(X_{\tau^*}\vert Z^*_{[0,\tau^*]})= \Lambda((-\iy,+\iy),\cdot)=\mu$
does not depend on $Z^*_{[0,\tau^*]}$.
It follows that 
$\tau^*$  is strong because it is measurable with respect to $Z^*_{[0,\tau^*]}$.
Finally it is a strong stationary time for $X$, since from the above identity,
$\cL(X_{\tau^*})=\mu$.\par\sm
Up to the construction of the intertwining, these few standard arguments 
provide the direct implication in Theorem \ref{t1}:
\begin{pro}\label{Itau}
If $I<+\iy$, then the random time $\tau^*$ defined in (\ref{taustar}) is a.s.\ finite and by consequence it is a strong stationary time for the positive recurrent diffusion $X$.
\end{pro}
\par
\begin{rem}
In Cheng and Mao \cite{Cheng_passage},  a strong stationary time is also obtained duality, up to the construction of the intertwining coupling.
In the situation where the process starts from one of the reflecting boundaries, the dual process is easier to deduce, because it is itself a one-dimensional diffusion (see also \cite{MR2530103}, which deals with strong quasi-stationary times for finite
birth and death process, but whose formalism is adapted to treat diffusion processes starting from the boundary).
\end{rem}
It opens the way to a quantitative study of the convergence to equilibrium for $X$
in the separation sense. Let us recall that the separation discrepancy $\mathfrak{s}(\nu,\mu)$ between two probability
measures $\nu$ and $\mu$ defined on the same state space $E$ is given by
\bq
\mathfrak{s}(\nu,\mu)&\df& \sup_{x\in E}1-\frac{d\nu}{d\mu}(x)\eq
where $\frac{d\nu}{d\mu}$ is the Radon-Nikodym derivative of the absolutely continuous part of $\nu$ with respect to $\mu$.
Strictly speaking, the separation discrepancy is not  a distance because it is not symmetric in its arguments.
The computations of Aldous and Diaconis \cite{MR876954}  show that 
for any strong stationary time $\tau$ for $X$, we have
\bqn{sharp}
\fo t\geq 0,\qquad \mathfrak{s}(\cL(X_t),\mu)&\leq & \PP[\tau >  t]\eqn
\par
Thus Proposition \ref{Itau} enables to get upper bounds on the speed of convergence of $X$
toward its equilibrium $\mu$ in the separation sense, by studying the speed of absorption 
at $(-\iy,+\iy)$ of $Z^*$. 
The inequalities (\ref{sharp}) may be equalities for all times $t\geq 0$ and such times $\tau$ are then stochastically minimal among all
strong stationary times. They are called sharp stationary times in 
Diaconis and Fill \cite{MR1071805} (in the finite setting).
The proof of the converse implication in Theorem \ref{t1} will rely on the fact that for initial distributions of $X$
of the form $\Lambda((-\iy,x),\cdot)$ and $\Lambda((x,+\iy),\cdot)$, with $x\in\RR$, the random time 
$\tau^*$ defined in (\ref{taustar}) is indeed a sharp stationary time.
\par
When is this technique working? It is convenient to consider the case of Langevin diffusions,
where $a\equiv 1$ and $b=-U'$, where $U\st\RR\ri\RR$ is a smooth potential.
In dimension 1 and up to shrinking the state space $\RR$ to an open interval (through a smooth transformation), it is not really a restriction.
The invariant measure $\mu$ admits  then a density proportional to $\exp(-U)$.
An application of Fubini's formula shows that the condition $I<+\iy$ writes down
\bqn{IU}
\max\lt(\int_{-\iy}^0\mu((-\iy,x))\frac1{\mu(x)}\, dx,\int_0^{+\iy}\mu((x,+\iy))\frac1{\mu(x)}\, dx\rt)&<& +\iy\eqn
\begin{rem}
The l.h.s.\ of (\ref{IU}) is bounded below by
\bq
\max\lt(\sup_{y\leq 0}\int_y^0\frac1{\mu(x)}\, dx\mu((-\iy,y)),\, \sup_{y\geq 0}\int_0^y\frac1{\mu(x)}\, dx\mu((y,+\iy))\rt)
\eq
and if $0$ was chosen to be the median of $\mu$ (up to a translation there is not lack of generality in this choice), 
the previous quantity is the inverse of the spectral gap of $L$ in $\LL^2(\mu)$ up to a factor 4 (see e.g.\ 
Bobkov and G{\"o}tze \cite{MR1682772}). So at least for Langevin diffusions, the existence of a strong stationary time,
whatever the initial distribution, implies a positive spectral gap. As it will appear below and as it can be expected from Remark
\ref{Mao}, this is far from being a sufficient
condition.
\end{rem}\par
For instance, if for $\lve x\rve$ large enough we have $U(x)=\lve x\rve^{\alpha}$, with $\alpha >0$, then Condition (\ref{IU})
is satisfied if and only if $\alpha>2$ (whereas the existence of a spectral gap is equivalent to $\alpha\geq 1$).
In particular, the important case of the Ornstein-Uhlenbeck process is not covered. Does it mean that the previous approach is
useless in this situation? Indeed, it is possible to get around this difficulty by considering strong times $\tau$ where the 
distribution of $X_\tau$ is close to the invariant probability $\mu$. Put in practice in Section 5,
this technique will enable us to recover good quantitative bounds on the convergence of the Ornstein-Uhlenbeck
process toward the Gaussian distribution in the total variation sense.
\par\sm
Let us just give a glimpse of why it could interesting to investigate the multidimensional situation.
Let $X$ be a hypoelliptic diffusion taking values in a smooth manifold $M$ of dimension (strictly) larger than 1. 
Assume that
it is possible to construct a process $Z^*$ taking values in the set $E^*$ of singletons and non-empty open subsets of $M$ and which is intertwined
with $X$ through the Markov kernel $\Lambda$ from $E^*$ to $M$ given  by
\bq
\fo z\in E^*,\qquad \Lambda(z, \,\cdot \,)&\df& 
\lt\{\begin{array}{ll}
\delta_{x}( \,\cdot \,)&\hbox{, if $z=\{x\}$}\\
\noalign{\vskip 2mm}
\frac{\lambda( \,\cdot \,\cap z)}{\lambda(z)}
&\hbox{, if $z$ is a non-empty subset of $M$}
\end{array}\rt.
\eq
where $\lambda$ is a nice $\sigma$-finite measure on $M$ giving positive weights to 
all non-empty open subsets (for instance the invariant measure for $X$, but it could also be a more tractable measure).
Then we would have at our disposal the following representation of the time marginal laws of $X$ for all $t\geq 0$,
\bq
\fo x\in M,\qquad\cL(X_t)(dx)&=&\int\Lambda(z,dx)\, \cL(Z^*_t)(dz)\eq
from which absolute continuity and regularity properties can be deduced.
\\
It would be instructive to begin with a simple instance of $X$ satisfying H\"ormander's conditions
and to see which features could be deduced for corresponding processes $Z^*$, especially in small times.
Entrance boundary properties of singletons analogous to that presented in Proposition \ref{entrance}
would be particularly desirable.
\par\me
The paper is constructed on the following plan. In the next section we investigate the dual process $Z^*$, making a link with the
square Bessel process of dimension 3 and 
we prove 
Proposition~\ref{entrance}. Explosion times and Proposition \ref{Itau} 
are the subject of Section 3.
Section 4  ends the proof of Theorem \ref{t1}, 
providing the missing details about the coupling of $X$ with $Z^*$ and 
showing the converse implication. 
 The last section and an appendix  are devoted to the counter-example
 of the benchmark Ornstein-Uhlenbeck process, 
  giving us the opportunity
 to see why it is interesting to consider more general
strong times than strong stationary times.

\section{Description of the dual process}\label{dotdp}

We study here the solutions of the stochastic differential equations associated with the generator $L^*$ given by (\ref{Lstar}), (\ref{Lstar2}) and (\ref{Lstar3}).\par\me
We begin by verifying the assertion made in the introduction about the relation between $L^*$ and $\wi L$ defined in (\ref{wiL}).
\begin{lem}\label{evolset}
Let $h$ be the function introduced in (\ref{h}). On $\rE$ we have
$\wi L[h]=0$ and for any $F\in\cC^\iy(\rE)$,
\bq
\fo z\in \rE,\qquad
L^*[F](z)&=&\frac1h\wi L[hF](z)\eq
These properties extend to $\RR\times\{+\iy\}$ and $\{-\iy\}\times\RR$, up to the natural modifications.
\end{lem}
\proof
For $(x,y)\in \rE$, we have
\bq
\pa_xh(x,y)\ =\ -\mu(x)&\hbox{ and }& \pa_yh(x,y)\ =\ \mu(y)\eq
so that
\bq
(\sqrt{a(y)}\pa_y-\sqrt{a(x)}\pa_x)^2h(x,y)
&=&(\sqrt{a(y)}\pa_y-\sqrt{a(x)}\pa_x) (\sqrt{a(y)}\mu(y)+\sqrt{a(x)}\mu(x))\\
&=&\sqrt{a(y)}\pa_y( \sqrt{a(y)}\mu(y))-\sqrt{a(x)}\pa_x (\sqrt{a(x)}\mu(x))\eq
Taking into account that
\bq
\mu'\ =\ \lt(c'-\frac{a'}{a}\rt)\mu\ =\ \frac{b-a'}{a}\mu\eq
we get that
\bq
(\sqrt{a(y)}\pa_y-\sqrt{a(x)}\pa_x)^2h(x,y)&=&
\lt(b(y)-\frac{a'(y)}{2}\rt)\mu(y)-\lt(b(x)-\frac{a'(x)}{2}\rt)\mu(x)\\
&=&-(a'(x)/2-b(x))\pa_xh(x,y)-(a'(y)/2-b(y))\pa_yh(x,y)\eq
namely $\wi L[h]=0$.\\
In the same way one shows that $\wi L[h]=0$ on  $(\RR\times\{+\iy\})\sqcup(\{-\iy\}\times\RR)$.
\par\sm
By definition of $\wi \Gamma$, we observe that 
for any $F\in\cC^\iy(\rE)$ and any $ z\in \rE$,
\bq
\frac1h\wi L[hF](z)&=&\frac1h(h\wi L[F]+F\wi L[h]+\wi \Gamma[h,F])\\
&=&\wi L[F]+\frac1h \wi \Gamma[h,F]\eq
A direct computation shows that for any $F,G\in\cC^\iy(\rE)$ and any $ z=(x,y)\in \rE$,
\bq
 \wi \Gamma[G,F]&=&2(\sqrt{a(y)}\pa_y-\sqrt{a(x)}\pa_x)G(x,y)
 (\sqrt{a(y)}\pa_y-\sqrt{a(x)}\pa_x)F(x,y)\eq
Applying this formula with $G=h$, we obtain that
$L^*[F]=\frac1hL[hF]$, as announced.\\
Again these considerations extend without difficulty to $(\RR\times\{+\iy\})\sqcup(\{-\iy\}\times\RR)$.
\wwtbp
\begin{rem}\label{dddagger}
Similar computations are valid for the generators given by (\ref{dagger}) and (\ref{ddagger}).
Indeed, they are respectively the Doob transforms through $h$ of the generators defined by
\bqn{checkL}
\check L&\df&(\sqrt{a(y)}\pa_y+\sqrt{a(x)}\pa_x)^2+(a'(x)/2-b(x))\pa_x+(a'(y)/2-b(y))\pa_y\eqn
and
\bqn{breveL}
\breve L&\df&a(y)\pa^2_y+a(x)\pa^2_y+(a'(x)/2-b(x))\pa_x+(a'(y)/2-b(y))\pa_y
\eqn
(on $\rE$ and their natural extensions on $\RR\times\{+\iy\}$ and $\{-\iy\}\times\RR$).
\\
Essentially relying on the fact that $\pa_x\pa_y h=0$, one deduces $\check L[h]=0$ from
$\wi L[h]=0$ and $\breve L[h]=0$ from $\breve L=(\wi L+\check L)/2$.
\end{rem}
\par
Note that the generator $L^*$ described in (\ref{Lstar}) expands into
\bq
 L^*&=&
a(x)\pa_x^2+a(y)\pa_y^2-\sqrt{a(x)}\sqrt{a(y)}\pa_x\pa_y+(a'(x)-b(x))\pa_x+(a'(y)-b(y))\pa_y\\&&+2\frac{\sqrt{a(x)}\mu(x)+\sqrt{a(y)}\mu(y)}{\mu([x,y])}(\sqrt{a(y)}\pa_y-\sqrt{a(x)}\pa_x)
\eq
It follows that
on $\rE$, the stochastic differential equation for $Z^*=(X^*,Y^*)$ associated to (\ref{Lstar})
writes down 
\bq
dX^*_t&=&\lt(a'(X^*_t)-b(X^*_t)-2\frac{\sqrt{a(X^*_t)}\mu(X^*_t)+\sqrt{a(Y^*_t)}\mu(Y^*_t)}{\mu([X^*_t,Y^*_t])}\sqrt{a(X^*_t)}\rt)dt -\sqrt{2a(X^*_t)}\,dB_t\\
dY^*_t&=&\lt(a'(Y^*_t)-b(Y^*_t)+2\frac{\sqrt{a(X^*_t)}\mu(X^*_t)+\sqrt{a(Y^*_t)}\mu(Y^*_t)}{\mu([X^*_t,Y^*_t])}\sqrt{a(Y^*_t)}\rt)dt +\sqrt{2a(Y^*_t)}\,dB_t
\eq
where $B=(B_t)_{t\geq 0}$ is a standard (one dimensional) Brownian motion. Starting from an initial condition in $\rE$, the regularity of the coefficients
and standard results (see for instance the book of Ikeda and Watanabe \cite{MR1011252})
show that there are existence and uniqueness of the solution $Z^*$ up to the explosion time $\tau^{\dagger}$ (a.s.\ with respect to
$B$). This stopping time for $Z^*$ is defined by
\bqn{taudagger}
\tau^\dagger&\df&\min( \tau_1,\tau_2,\tau_3)\eqn
where
\bq
\tau_1&\df&\lim_{r\ri +\iy}\inf\{t\geq 0\st X^*_t< -r\}\\
\tau_2&\df&\lim_{r\ri +\iy}\inf\{t\geq 0\st Y^*_t > r\}\\
\tau_3&\df&\lim_{r\ri +\iy}\inf\{t> 0\st Y^*_t-X^*_t<1/r\}\eq
Of course, we have $\tau^\dagger\leq \tau^*$, where $\tau^*$ is defined in (\ref{taustar}).
The next result shows that $\tau_3$ plays no role.
\begin{lem}\label{mart}
Let $Z^*$ start from an initial condition in $\rE$. Then a.s.\ $h(Z^*_t)$ converges as $t$ goes to $\tau^\dagger$
toward a positive quantity. In particular $\tau^\dagger=\min(\tau_1,\tau_2)$ and $Z^*$ can exit $\rE$
only through $(\RR\times\{+\iy\})\sqcup(\{-\iy\}\times\RR)$.
\end{lem}
\proof
According to Lemma \ref{evolset}, we have on $\rE$, $L^*[1/h]=\wi L(\un)/h=0$, where $\un$ is the function always taking the value
1 on $\rE$. It follows that the process $M=(M_t)_{t\geq 0}$ defined by
\bq
\fo t \geq 0,\qquad 
M_t&\df& \frac{1}{h(Z^*_{\tau^\dagger\wedge t})}\eq
is a local martingale. Since it is furthermore positive, it must converge as $t$ goes to infinity.
The announced results follow.\wwtbp
\par
We can now obtain the equivalent of Proposition \ref{entrance} but for initial conditions in $\rE$.
\begin{pro}\label{facile}
For any $z_0\in\rE$, there is a continuous Markov process $Z^*\df(Z^*_t)_{t\geq 0}$ starting from $z_0$ and whose generator is $L^*$.
The law of this process is unique if we impose that after the possibly finite time
$
\tau^*$, defined as in (\ref{taustar}), $Z^*$ stays at  position $(-\iy,+\iy)$ . Furthermore for all $t\geq 0$, $Z_t\in E^*\setminus D^*$.
\end{pro}
\proof
According to the previous arguments, we already have the existence and uniqueness of $Z^*$ up
to the time $\tau^\dagger$. If $\tau^\dagger=+\iy$, the construction is over.
If $\tau^\dagger<+\iy$, we deduce from Lemma \ref{mart} that either 
$\tau_1=\tau^\dagger<+\iy$, or $\tau_2=\tau^\dagger<+\iy$. We only consider the first case, the second can be treated in the same way.
By the required continuity of the trajectories, we must have $Z^*_{\tau^\dagger}=(X^*_{\tau_1},+\iy)$, where $X^*_{\tau_1}\in [-\iy,+\iy)$.
We first consider the case where $X^*_{\tau_1}\not=-\iy$. By the assumption on the form of $L^*$ on $\RR\times\{+\iy\}$, $Z^*$ must stay there after time $\tau^\dagger$.
Let us denote for any $t\geq 0$, $X^{\dagger}_t\df X^*_{\tau_1+t}$. The process $X^{\dagger}$
must be (and is constructed as) a solution of the one-dimensional stochastic differential equation
\bq
dX^\dagger_t&=&\lt(a'(X^\dagger_t)-b(X^\dagger_t)-2\frac{a(X^\dagger_t)\mu(X^\dagger_t)}{\mu([X^\dagger_t,+\iy))}\rt)dt -\sqrt{2a(X^\dagger_t)}\,dB_t\eq
(where $B=(B_t)_{t\geq 0}$ is a standard  Brownian motion), starting from $X^*_{\tau_1}$. 
Due to the regularity of the coefficients, there is no difficulty to get existence and uniqueness of  the solution up to the time
\bq
\bar \tau&\df&\lim_{r\ri +\iy}\inf\{t\geq 0\st \lve X^\dagger_t\rve > r\}\eq
As in the proof of Lemma \ref{mart}, the process 
$M^\dagger=(M^\dagger_t)_{t\geq 0}$ defined by
\bq
\fo t \geq 0,\qquad 
M^\dagger_t&\df& \frac{1}{h((X^\dagger_{\bar\tau\wedge t},+\iy))}\eq
is a positive local martingale. From its convergence we deduce that
\bq
\lim_{t\ri \bar\tau-}X^{\dagger}_t&=&-\iy\eq
and it follows that $\tau^{\dagger}+\bar\tau=\tau^*$. Note that this identity is trivial if $X^*_{\tau_1}=-\iy$.
The analogue result  is satisfied in the situation $\tau_2=\tau^\dagger$.
Thus the law of $Z^*_{[0,\tau^*)}$ is uniquely determined and since we impose that
$Z^*_t=(-\iy,+\iy)$ for $t\geq \tau^*$ (by continuity  for $t=\tau^*$), the same is true for $Z^*$.
The fact that $Z_t\in E^*\setminus D^*$ for all $t\geq 0$ is obvious from the previous martingale arguments.\wwtbp
\par
For $z_0\in\rE$, designate by $\PP_{z_0}$ the law on the set of trajectories $\cC(\RR_+,E^*)$ of $Z^*$ starting from $z_0$ and constructed as above.
One way to construct $\PP_{z_0}$ for $z_0=(x_0,x_0)\in\DD^*$, is to consider for $\epsilon,\epsilon' >0$,
$\PP_{x_0-\epsilon,x_0+\epsilon'}$ and to let $\epsilon,\epsilon'$ go to zero.
To make clearer the convergence, we will consider a transformation of  $\cC(\RR_+,E^*)$
so that all the difficulties are encapsulated into a square Bessel process of dimension 3.\par
Here is how it appears: under $\PP_{z_0}$ for some $z_0\in \rE$, consider
\bqn{varsigma}
\varsigma&\df&2\int_0^{\tau^*}(\sqrt{a(X^*_s)}\mu(X^*_s)+\sqrt{a(Y^*_s)}\mu(Y^*_s))^2\,ds\ \in(0,+\iy]\eqn
(with the convention $\sqrt{a(\pm\iy)}\mu(\pm\iy)=0$),
and
 the time change $(\theta_t)_{t\in[0,\varsigma]}$
defined by
\bqn{theta}
\fo t\in[0,\varsigma],\qquad2 \int_0^{\theta_t}(\sqrt{a(X^*_s)}\mu(X^*_s)+\sqrt{a(Y^*_s)}\mu(Y^*_s))^2 \,ds&=&t\eqn
We are interested in the process $R\df(R_t)_{t\geq 0}$ given by
\bqn{R}
\fo 
 t\geq 0,\qquad R_t&\df& h(Z^*_{\theta_{t\wedge \varsigma}})\eqn
 \begin{pro}\label{Bessel}
 Under $\PP_{z_0}$ with $z_0\in \rE$, $R$ has the law of a 
square Bessel process of dimension 3 starting from $h(z_0)\in(0,1)$
and stopped at 1. In particular $\varsigma$ is distributed as the first reaching time of 1 for this process.
\end{pro}
\proof
We begin by computing $L^*[h]$: in view of Lemma \ref{evolset} we have on $\rE$,
\bq
L^*[h]&=& \frac1h\wi L[h^2]\\
&=& \frac1h(2h\wi L[h]+\wi\Gamma[h,h])\\
&=&\frac1h\wi\Gamma[h,h]\eq
Taking into account the stochastic differential equations satisfied by the coordinates $X^*$ and $Y^*$ of $Z^*$,
Itô's formula give us
\bq
dh(Z_t^*)&=&\frac1{h(Z^*_t)}\wi\Gamma[h,h](Z^*_t)\,dt +(\sqrt{2a(X_t^*)}\mu(X_t^*)+\sqrt{2a(Y_t^*)}\mu(Y_t^*))\,dB_t\eq
In Lemma \ref{evolset} we have already seen that
\bq
\fo z=(x,y)\in\rE,\qquad
\wi \Gamma[h,h](z)&=&2(\sqrt{a(x)}\mu(x)+\sqrt{a(y)}\mu(y))^2\eq
Classical stochastic time change calculus (cf.\ for instance Chapter 5 of the book \cite{MR1725357} of Revuz and Yor)
then shows that the process $R$ satisfies for $t<\varsigma$:
\bq
dR_t&=&\frac{d\theta_t}{dt}\frac1{h(Z^*_{\theta_t})}\wi\Gamma[h,h](Z^*_{\theta_t})\,dt +\sqrt{\frac{d\theta_t}{dt}}\sqrt{\wi\Gamma[h,h](Z^*_{\theta_t})}\,dW_t\eq
where $W=(W_t)_{t\geq 0}$ is a standard  Brownian motion. From the definition of the time change $(\theta_t)_{t\in[0,\varsigma)}$,
we have
\bq
\fo t\in[0,\varsigma),\qquad \frac{d\theta_t}{dt}&=&\frac1{\wi\Gamma[h,h](Z^*_{\theta_t})}\eq
so we end up with
\bq
dR_t&=&\frac1{R_t}\, dt+dW_t\eq
One recognizes the stochastic differential equation characterizing the square Bessel process of dimension 3 (see e.g.\ Chapter 11 of the book \cite{MR1725357} of Revuz and Yor).
Since $Z^*$ is stopped when it reaches $(-\iy,+\iy)$, namely when $h(Z^*)$ hits 1, $R$ is stopped when it reaches 1,
which ends to show the assertions of the proposition.\wwtbp
\par 
Here is a first consequence of the previous result:
\begin{cor}\label{convZ}
We have almost surely,
\bq
\lim_{t\ri \tau^*-}X^*_{t}&=&-\iy\\
\lim_{t\ri \tau^*-}Y^*_{t}&=&+\iy\eq
\end{cor}
\proof
From (\ref{varsigma}) and (\ref{theta}), we get that 
as $t$ converges to $\varsigma-$, $\theta_t$ converges to $\tau^*-$.
It follows that
\bq
\lim_{t\ri \tau^*-}h(Z_t^*)\ =\ \lim_{t\ri \varsigma-} R_t\ =\ 1\eq
Recalling the definition of $h$ given in (\ref{h}), this is possible
if and only if the limits described in the above corollary take place.\wwtbp\par
The idea behind the proof of Proposition \ref{entrance} is that there is no difficulty to let a square Bessel process of dimension 3
start from 0. But to proceed rigorously, we need to consider some transformations of the martingale problem
associated to $L^*$.
First we remark that it is sufficient to show that 
for any $x_0\in\RR$, there is a continuous Markov process $(Z^*_t)_{t\in[0,\tau^\dagger)}$,  where $\tau^\dagger$ is defined as in (\ref{taudagger}), starting from $(x_0,x_0)$, living in $E^*\setminus D^*$ for $t\in (0,\tau^\dagger)$,
and 
whose generator is $L^*$.
Furthermore, we will check that the law of this process is unique.
Indeed, the proof of Proposition \ref{facile} could next be used again to uniquely extend $(Z^*_t)_{t\in[0,\tau^\dagger)}$
into $(Z^*_t)_{t\geq 0}$.
This observation brings us back to the martingale problem associated to the initial condition $(x_0,x_0)$
and to the restriction of the generator $L^*$ to $\rE$.
But we begin by replacing $(x_0,x_0)$ by $(x_0-\epsilon,x_0+\epsilon')$, with $\epsilon,\epsilon'>0$, and we consider the time change described in (\ref{theta}).
This amounts to replace the generator $L^*$ by $\wit L\df (1/\wi \Gamma(h,h))L^*$.
Or equivalently to apply the following transformation to the trajectories
\bq
(Z^*_t)_{t\in[0,\tau^\dagger)}&\mapsto& (\wit Z_t)_{t\in[0,\varsigma)}\df (Z^*_{\theta_t})_{t\in[0,\varsigma)}\eq
where 
\bq
\varsigma&\df&2\int_0^{\tau^\dagger}(\sqrt{a(X^*_s)}\mu(X^*_s)+\sqrt{a(Y^*_s)}\mu(Y^*_s))^2\,ds\ \in(0,+\iy]\eq
and $(\theta_t)_{t\in[0,\varsigma)}$ is defined as in  (\ref{theta}).
The reverse mapping is given by
\bq
(\wit Z_t)_{t\in[0,\varsigma)}&\mapsto&(Z^*_t)_{t\in[0,\tau^\dagger)} \df (\wit Z_{\vartheta_t})_{t\in[0,\tau^\dagger)}\eq
where
\bq
\tau^{\dagger}&\df& 
\frac12\int_0^{\varsigma}(\sqrt{a(\wit X_s)}\mu(\wit X_s)+\sqrt{a(\wit Y_s)}\mu(\wit Y_s))^{-2}\,ds\ \in(0,+\iy]\eq
and
\bq
\fo t\in[0,\tau^{\dagger}),\qquad\frac12 \int_0^{\vartheta_t}(\sqrt{a(\wit X_s)}\mu(\wit X_s)+\sqrt{a(\wit Y_s)}\mu(\wit Y_s))^{-2} \,ds&=&t\eq
The
 stochastic differential equation for $(\wit Z_t)_{t\in[0,\varsigma)}=(\wit X_t,\wit Y_t)_{t\in[0,\varsigma)}$ associated to $\wit L$ on $\rE$ is given by
\bq
\fo t\in[0,\varsigma),\qquad\lt\{\begin{array}{rcl}
d\wit X_t&=&b_1(\wit X_t,\wit Y_t)\,dt +\sigma_1(\wit X_t,\wit Y_t)\,dB_t\\
d\wit Y_t&=&b_2(\wit X_t,\wit Y_t)\,dt +\sigma_2(\wit X_t,\wit Y_t)\,dB_t\end{array}\rt.
\eq
where for any $(x,y)\in \rE$,
\bq
b_1(x,y)&\df& \frac{a'(x)-2b(x)}{4(\sqrt{a(y)}\mu(y)+\sqrt{a(x)}\mu(x))^2)}
-\frac{\sqrt{a(x)}}{(\sqrt{a(y)}\mu(y)+\sqrt{a(x)}\mu(x))\mu([x,y])}\\
b_2(x,y)&\df& \frac{a'(y)-2b(y)}{4(\sqrt{a(y)}\mu(y)+\sqrt{a(x)}\mu(x))^2)}
+\frac{\sqrt{a(y)}}{(\sqrt{a(y)}\mu(y)+\sqrt{a(x)}\mu(x))\mu([x,y])}\\
\sigma_1(x,y)&\df& -
\frac{\sqrt{a(x)}}{\sqrt{a(y)}\mu(y)+\sqrt{a(x)}\mu(x)}\\
\sigma_2(x,y)&\df& 
\frac{\sqrt{a(y)}}{\sqrt{a(y)}\mu(y)+\sqrt{a(x)}\mu(x)}\eq
\par
Finally we consider the transformation $\Psi$ of the state space $\rE\sqcup D^*$ given by
\bq
\rE\sqcup D^*\ni (x,y)&\mapsto & (h(x,y),s(x,y))\eq
where $s(x,y)$ is the middle point of $[x,y]$ when $\RR$ is endowed with the Riemannian structure 
for which $a\pa^2$ is the Laplace-Beltrami operator. More prosaically, $s(x,y)$ is defined as the unique point
in $[x,y]$ such that
\bq
\int_x^{s(x,y)}\frac{1}{\sqrt{a(u)}}\,du&=&\int_{s(x,y)}^y\frac{1}{\sqrt{a(u)}}\,du\eq
Its main interest is that 
\bq
\fo (x,y)\in\rE,\qquad \sqrt{a(y)}\pa_y s(x,y)-\sqrt{a(x)}\pa_x s(x,y)&=&0\eq
because
\bqn{sxy}
\pa_xs(x,y)\ =\ \frac12\sqrt{\frac{a(s(x,y))}{a(x)}}&\hbox{ and }&\pa_ys(x,y)\ =\ \frac12\sqrt{\frac{a(s(x,y))}{a(y)}}\eqn
It is not difficult to see that $\Psi$ is a smooth diffeomorphism from $\rE\sqcup D^*$ to its image.
Denote  $(R_t,S_t)_{t\in[0,\varsigma)}\df (\Psi (\wit Z_t))_{t\in[0,\varsigma)}$.
From Proposition \ref{Bessel} and (\ref{sxy}) we deduce that
the stochastic differential equation satisfied by $\wit Z_{t\in[0,\varsigma)}$ is transformed into
\bqn{RS}
\fo t\in[0,\varsigma),\qquad\lt\{\begin{array}{rcl}
dR_t&=&\frac1{R_t}\, dt+dW_t\\
dS_t&=&\beta(R_t,S_t)\,dt\end{array}\rt.
\eqn
where $W=(W_t)_{t\geq 0}$ is a standard  Brownian motion
 and where the mapping $\beta$ is defined on $ \Psi(\rE\sqcup D^*)$ by
 \bq
 \fo(x,y)\in\rE\sqcup D^*,\qquad \beta(\Psi(x,y))&\df&
  \frac{a'(x)-2b(x)}{8(\sqrt{a(y)}\mu(y)+\sqrt{a(x)}\mu(x))^2}\sqrt{\frac{a(s(x,y))}{a(x)}}\\&&+
  \frac{a'(y)-2b(y)}{8(\sqrt{a(y)}\mu(y)+\sqrt{a(x)}\mu(x))^2}\sqrt{\frac{a(s(x,y))}{a(y)}}
 \eq
It is clear that this function is smooth on its domain.
So the resolution of (\ref{RS}) is quite obvious. The initial condition is $(R_0,S_0)
=\Psi((x_0-\epsilon,x_0+\epsilon'))$.
Then one solves the autonomous stochastic differential equation satisfied by $R\df(R_t)_{t\geq 0}$.
The solution $R$ is defined for all $t\geq 0$, and as it was more precisely seen in Proposition \ref{Bessel}, it is a square Bessel process of dimension 3 starting from
$h(x_0-\epsilon,x_0+\epsilon')>0$. The trajectory $R$ being constructed,
it remains to investigate the ordinary differential equation $\frac{dS_t}{dt}=\beta(R_t,S_t)$, starting from $S_0$.
Since $\beta$ is smooth, it gives us a solution, up to the possible explosion time $\varsigma$ when $(R_{\varsigma-},S_{\varsigma-})$
reaches the boundary of $ \Psi(\rE\sqcup D^*)$.
From the form of $\Psi$, the time $\varsigma$ is necessarily the first time when
either
$\wit X$ explodes to $-\iy$ or $\wit Y$ explodes to $+\iy$, where
$(\wit X_t,\wit Y_t)\df\Psi^{-1}(R_t,S_t)$ for $t\in[0,\varsigma)$, as wanted.\par
These observations are also valid if $R_0=0$ and enable to construct 
$\PP_{(x_0,x_0)}$
by reversing the previous transformations, starting from the initial condition $(R_0,S_0)=(0,x_0)$.
It is also seen to be the limit of $\PP_{(x_0-\epsilon,x_0+\epsilon)}$ as $\epsilon,\epsilon'>0$ converge to zero.
In the last sentence, the weak convergence of the probability measures is with respect to the uniform convergence of the trajectories over compact time intervals,
when the state space $E^*$ is endowed with a bounded distance compatible with its natural topology (inherited from that of the compact set 
$[-\iy,+\iy]^2$).
This continuity property  and the requirements made on $\PP_{(x_0,x_0)}$ in Proposition \ref{entrance}
enable  us to be convinced
of its uniqueness.
Indeed, consider $\PP$ another probability on the trajectories $\cC(\RR_+,E^*)$ satisfying the same properties.
On $\cC(\RR_+,E^*)$ consider the  natural time-shift maps $\Theta_t$,
for $t\geq 0$:  if $(Z^*_s)_{s\geq 0}$ stands for the canonical coordinate process, we have
$Z^*_s(\Theta_t)=Z^*_{t+s}$ for all $t,s\geq 0$.
Let $F$ be a bounded and continuous functional on $\cC(\RR_+,E^*)$.
By the Markov property we must have that for any $t>0$,
\bq
\EE[F(\Theta_t)]&=&\int \EE_{z}[F]\,m_t(dz)\eq
where $m_t$ is the law of $Z_t^*$ under $\PP$.
By the requirements that $Z^*_t$ belongs to $\EE^*\setminus D^*$ a.s.\ under $\PP$ and that $\PP$ is a solution
to the martingale problem associated to $L^*$,
the expectations $\EE_z$ in the r.h.s.\ are relatively to the laws constructed in Proposition \ref{facile}.
By the continuity of the trajectories and the assumption that $Z_0^*=(x_0,x_0)$ under $\PP$,
$m_t$ converges weakly to the Dirac mass at $(x_0,x_0)$ as $t$ goes to $0_+$.
Thus $\lim_{t\ri 0_+}\int \EE_{z}[F]\,m_t(dz)= \EE_{(x_0,x_0)}[F]$ by the continuity of $z\mapsto \PP_z$ at $(x_0,x_0)$.
On the other hand, by the dominated convergence theorem, $\lim_{t\ri 0_+}\EE[F(\Theta_t)]=\EE[F]$.
Thus $\EE[F]=\EE_{(x_0,x_0)}[F]$ for all bounded and continuous functional $F$ on $\cC(\RR_+,E^*)$.
This is sufficient to insure that $\PP=\PP_{(x_0,x_0)}$ and ends the proof of Proposition \ref{entrance}.
\par
\begin{rem}\label{alpha}
Proposition \ref{facile} and its proof are also valid for the generators defined in (\ref{dagger}) and (\ref{ddagger}) and more generally for
the generators
$L_{\alpha}^*\df(1-\alpha) L^*+ \alpha \check L^*$, where  $\alpha\in[0,1]$. But 
Proposition \ref{entrance} is not true for $\check L^{*}$: 
as it was mentioned in the introduction, due to the regularity of the coefficients of $\check L^{*}$,
the unique solutions $\check \PP_z$  for the corresponding martingale problem
can be directly constructed for all the initial conditions $z\in E^*$ and the mapping $z\mapsto \PP_z$ is continuous.
Unfortunately, starting from $z\in D^*$, the process cannot escape from $D^*$, except by possibly exploding at one of its two ends
(Lemma \ref{mart} is not helping to prevent this event: $h(Z^*_t)$ remains null). Indeed in this degenerate situation one may have to
add the two absorbing points $(-\iy,-\iy)$ and $(+\iy,+\iy)$ to the state space $E^*$.\par
This problem is not encountered 
by the generators $L_{\alpha}^*$,
for $\alpha\in [0,1)$, to which the above considerations (corresponding to the case $\alpha=0$) can be extended.
Let us put  a corresponding index $\alpha$ to all the objects we have considered so far when $L^*$ is replaced by $L_\alpha^*$. For instance
we introduce the generator $\wi L_{\alpha}\df (1-\alpha)\wi L+ \alpha \check L$
and we compute that its carré du champ $\wi\Gamma_\alpha$ satisfies
\bqn{wGa}
\nonumber\fo z=(x,y)\in\rE,\qquad
\wi \Gamma_\alpha[h,h](z)&=&2(\sqrt{a(y)}\mu(y)+\sqrt{a(x)}\mu(x))^2-8\alpha \sqrt{a(x)a(y)}\mu(x)\mu(y)\\
&=& \wi \Gamma[h,h](z)-8\alpha  \sqrt{a(x)a(y)}\mu(x)\mu(y)
\eqn
It leads us to replace (\ref{varsigma}), (\ref{theta}) and (\ref{R}) respectively by
 \bqn{varsigmaa}
\varsigma_\alpha&\df&\int_0^{\tau^*}\wi \Gamma_\alpha[h,h](Z^*_s)\,ds\ \in(0,+\iy]\eqn
\bq
\fo t\in[0,\varsigma_\alpha],\qquad \int_0^{\theta_{\alpha,t}}\wi \Gamma_\alpha[h,h](Z^*_s) \,ds&=&t\eq
and
\bq
\fo 
 t\geq 0,\qquad R_{\alpha,t}&\df& h(Z^*_{\theta_{\alpha,t\wedge \varsigma_\alpha}})\eq
The interest of the latter process is that under $\PP_{\alpha,z_0}$, it is again a square process of dimension 3 starting from $h(z_0)$ and stopped at 1. The proof is identical to that of Proposition \ref{Bessel}.\par
But from (\ref{wGa}), we get that for any $z\in E^*\setminus D^*$, the quantity $\wi \Gamma_\alpha[h,h](z)$ is non-increasing in $\alpha \in [0,1)$ (it is decreasing when $z\in \rE$).
It follows from (\ref{varsigmaa}), that for any fixed $z_0\in \rE\sqcup D^*$,  if $\alpha_1<\alpha_2\in [0,1)$, then the law of $\tau^*$ under $\PP_{\alpha_1,z_0}$ is strictly larger than the law of $\tau^*$ under $\PP_{\alpha_2,z_0}$, with respect to the usual stochastic ordering
of laws on $\RR_+\sqcup\{+\iy\}$. 
Hence among all the generators  $L_{\alpha}^*$
for $\alpha\in [0,1)$, $L^*=L_0^*$ leads to the dual process $Z^*$ to be the fastestly absorbed at $(-\iy,+\iy)$ and thus is the most adequate for our purpose of constructing relatively small stationary times for $L$.
In words, the mirror-symmetry coupling of the Brownian motions at the boundary of the evolving segment
is optimal and the identical coupling is the worst (being utterly useless for the evolving segments starting from a singleton).
\end{rem}

\section{Explosion times}\label{Secei}

Our main objective here is to prove the finiteness assertion of Proposition \ref{Itau}. 
The arguments are based on comparisons with some appropriate diffusions on half-lines.\par\sm
Consider $Z^*=(X^*,Y^*)$ the process described in Proposition \ref{entrance} (for some fixed $x_0\in \RR$)
and constructed in the previous section. We are interested in the (total) explosion time $\tau^*$ defined in (\ref{taustar})
and our main task is to show that it is almost surely finite if $I<+\iy$.
So let us consider the (partial) explosion times
\bq
\tau^-&\df& \inf\{t\geq 0\st X^*_t=-\iy\}\\
\tau^+&\df& \inf\{t\geq 0\st Y^*_t=+\iy\}
\eq
For our purpose it is sufficient to show the following result (recall that $I_-$ and $I_+$ were defined  just above the statement of Theorem \ref{t1}).
\begin{pro}\label{Itau+}
If $I_+<+\iy$, then $\tau^+$ is a.s.\ finite.
\end{pro}
Indeed, by symmetry it will follow that
if $I_-<+\iy$, then $\tau^-$ is a.s.\ finite, so that $\tau^*=\tau^-\vee\tau^+<+\iy$ a.s.\ if $I<+\iy$.
\par
The proof of Proposition \ref{Itau+} relies on the comparison of $Y^*$ with a diffusion $U\df(U_t)_{t\geq 0}$
taking values
in $\RR_+\sqcup\{+\iy\}$, reflected at 0, absorbed at $+\iy$ and whose generator on $(0,+\iy)$ is $a\pa^2-(b-a' +2ak')\pa$, where $k$
is the mapping
$\RR\ni x\mapsto \ln(\mu((-\iy,x]))$.
More precisely, we take for $U$ the solution of the stochastic differential equation
\bq
dU_t&=&\lt(a'(U_t)-b(U_t)+2a(U_t)k'(U_t)\rt)dt +\sqrt{2a(Y^*_t)}\,dB_t+dl_t(U)
\eq
up to the explosion time $\tau(U)=\inf\{t\geq 0\st U_t=+\iy\}$, 
where $(l_t(U))_{t\geq 0}$ is the local time of $U$ at 0 and
where $B=(B_t)_{t\geq 0}$ is the same standard  Brownian motion as the one driving the s.d.e.\ satisfied by $Y^*$
\bq
dY^*_t&=&\lt(a'(Y^*_t)-b(Y^*_t)+2\frac{\sqrt{a(X^*_t)}\mu(X^*_t)+\sqrt{a(Y^*_t)}\mu(Y^*_t)}{\mu([X^*_t,Y^*_t])}\sqrt{a(Y^*_t)}\rt)dt +\sqrt{2a(Y^*_t)}\,dB_t
\eq
for $t\leq \tau^+$ (with the natural modification of the drift term if $X^*_t=-\iy$).
The interest is that the quantity $[\sqrt{a(X^*_t)}\mu(X^*_t)+\sqrt{a(Y^*_t)]\mu(Y^*_t)}\sqrt{a(Y^*_t)}/{\mu([X^*_t,Y^*_t])}-2a(Y^*_t)k'(Y^*_t)$
is non-negative and even positive for $0<t<\tau^-$. So if  $U$ and $Y^*$ started from the same initial condition $u_0\in (0,+\iy)$,
then $U$ stays below $Y^*$ up to the time
\bq
T&\df& \inf\{t\geq 0\st U_t=0\}\eq
and this is true whatever the behavior of $X^*$:
\begin{lem}\label{UY}
For all $t\in[0,T]$, we have $U_t\leq Y^*_t$.\end{lem}
As usual, this assertion has to be understood a.s., but not to burden the presentation, this is assumed to be implicit from now on.
Note also that after the time $T$, the local time $(l_t(U))_{t\geq 0}$ starts to play a role and $U$ can end being above $Y^*$.
\proof
This kind of comparison result is standard, see for instance Section 1 of Chapter 6
of the book of Ikeda and Watanabe \cite{MR1011252}.
Nevertheless, we find 
more illuminating to present a simple and direct proof
than to check their assumptions via localizing 
 arguments.\par
It is convenient to first transform $\RR\sqcup\{-\iy,+\iy\}$ via the mapping $A$ given by
\bq
\fo u\in \RR_+\sqcup\{+\iy\},\qquad A(u)&\df& \int_0^{u}\frac1{\sqrt{2a(v)}}\,dv\eq
Next we consider the processes $\wi U\df(\wi U_t)_{t\geq 0}\df
(A(U_t))_{t\geq 0}$ and $\wi Y\df(\wi Y_t)_{t\geq 0}\df
(A(Y^*_t))_{t\geq 0}$.
 Owing to Itô's formula, for $t\in[0,T\wedge\tau^+)$, they satisfies respectively the s.d.e.
\bq
d\wi U_t&=&f(\wi U_t)\,dt+dB_t\\
d\wi Y_t&=&(f(\wi Y_t)+S_t)\,dt+dB_t
 \eq
where 
\bq
\fo u\in\RR,\qquad f(u)&\df& \lt(a'(u)-b(u)-2a(u)k'(u)\rt)A'(u)+\frac12A''(u)\eq
and $S\df(S_t)_{t\geq 0}$ is the previsible process given by
\bq
\fo t\geq 0,\qquad S_t&\df& \sqrt{2}\frac{\sqrt{a(X^*_t)}\mu(X^*_t)+\sqrt{a(Y^*_t)}\mu(Y^*_t)}{\mu([X^*_t,Y^*_t])}-\sqrt{2a(Y^*_t)}k'(Y^*_t)\\
&\geq &  \sqrt{2}\frac{\sqrt{a(Y^*_t)}\mu(Y^*_t)}{\mu([X^*_t,Y^*_t])}-\sqrt{2a(Y^*_t)}\frac{\mu(Y_t^*)}{\mu((-\iy,Y_t^*])}\\
&=&\sqrt{2a(Y^*_t)}\frac{\mu(Y_t^*)\mu((-\iy,X_t^*))}{\mu([X^*_t,Y^*_t])\mu((-\iy,Y_t^*])}\\
&\geq &0
\eq
As already mentioned, what is important is this non-negativity of $S$.
Consider 
\bq
\sigma&\df&\inf\{t\in[0,T)\st \wi U_t>\wi Y_t\}\eq
with the usual convention that $\sigma\df+\iy$ if the set in the r.h.s.\ is empty.
We proceed by contradiction: assume that $\sigma<T$ (and in particular $\sigma$ is finite).
Necessarily we also have $\sigma<\tau^+$, because $\wi Y_t=A(+\iy)\geq A(U_t)=\wi U_t$ for all $t\geq \tau^+$.
By continuity $\wi U_\sigma=\wi Y_\sigma$ and
we consider two cases:\par
$\bullet$ If $\sigma<\tau^-$, then $S_{\sigma}>0$, thus there exists
$\epsilon >0$ such that for $s\in[\sigma,\sigma+\epsilon]$, 
$f(\wi Y_s)+S_s-f(\wi U_s)>0$. 
From the above s.d.e.\  we deduce that for all $\epsilon'\in(0,\epsilon]$,
\bq
\wi Y_{\sigma+\epsilon'}-\wi U_{\sigma+\epsilon'}
&=&\int_\sigma^{\sigma+\epsilon'}f(\wi Y_s)+S_s-f(\wi U_s)\, ds
\ >\ 0\eq
and this contradicts the definition of $\sigma$.
It follows that  $\wi U_t\leq \wi Y_t$ for all $t\in[0,T)$ and by continuity this is also true for $t=T$.
\par
$\bullet$ If $\sigma\geq \tau^-$: for  $t\geq \tau^-$, $S_t=0$, so $(\wi U_t)_{\tau^-\leq t\leq T\wedge\tau^+}$ and $(\wi Y_t)_{\tau^-\leq t\leq T\wedge\tau^+}$ follow the same
s.d.e.\ whose coefficients are regular. Since $\wi U_\sigma=\wi Y_\sigma$, the local uniqueness of the
solution of their s.d.e.\ implies that $\wi U$ and $\wi Y$ keep on being equal for some
time after $\sigma$ and this is again contradictory with the definition of $\sigma$.
\wwtbp
\par
The advantage of the process $U$ is that its explosion time $\tau(U)$ is well-understood, as we deduce from
Theorem 3.2 of Chapter 6 of the book of Ikeda and Watanabe \cite{MR1011252}
 the following criterion:
\begin{pro}\label{Uiy}
The explosion time $\tau(U)$ is finite almost surely if and only if $I_+<+\iy$.
\end{pro}
\proof
The most convenient way to exploit Section 3 of Chapter 6 of the book of Ikeda and Watanabe \cite{MR1011252}
seems to symmetrize $U$: consider the functions $\wit a$ and $\wit b$ defined by 
\bq
\fo x\in \RR,\qquad\lt\{\begin{array}{rcl}
\wit a(x)&\df& \lt\{
\begin{array}{ll}
a(x)&\hbox{, if $x\geq 0$}\\
a(-x)&\hbox{, if $x<0$}
\end{array}\rt.
\\
\wit b(x)&\df& \lt\{
\begin{array}{ll}
a'(x) -b(x)+2a(x)k'(x)&\hbox{, if $x\geq 0$}\\
-\wit b(-x)&\hbox{, if $x<0$}
\end{array}\rt.
\end{array}\rt.\eq
to which we associate the operator
\bq
\wit L&\df& \wit a\pa^2+\wit b\pa\eq
Since $\wit a$ is continuous and 
positive 
and $\wit b$
is measurable and locally bounded, we can use Theorem~3.3 of Chapter 4 of the book of Ikeda and Watanabe \cite{MR1011252}
and usual localization procedures to obtain, for  any given starting point $v\in\RR$, the existence and uniqueness 
of the solution $V\df(V_t)_{0\leq t\leq \tau(V)}$ of the s.d.e.\ associated to $v$ and $\wit L$:
\bq
\lt\{
\begin{array}{rcl}
V_0&=&v\\
dV_t&=&\wit b(V_t)dt +\sqrt{2\wit a(V^*_t)}\,dB_t
\end{array}\rt.
\eq
up to the explosion time $\tau(V)\df\inf\{t\geq 0\st\lim_{s\ri t-}\lve V_s\rve=+\iy\}$, 
 and
where $B=(B_t)_{t\geq 0}$ is a standard  Brownian motion.
\\
Tanaka's formula (e.g.\ Chapter 6 of the book \cite{MR1725357} of Revuz and Yor) enables to see that
$(\lve V_t\rve)_{0\leq t< \tau(V)}$ coincides in law with the process  $(U_t)_{0\leq t< \tau(U)}$
starting from $\lve v\rve$.
Formally, if
Theorem 3.2 (3) of Chapter 6 of the book of Ikeda and Watanabe \cite{MR1011252} is applied (take $c=0$ there),
we get that  the a.s.\ finiteness of $\tau(U)=\tau(V)$ (independently of the initial condition)
is equivalent to
\bqn{II}
\int_0^{+\iy}\exp\lt(-\int_0^x\frac{\wit b(y)}{\wit a(y)}\,dy\rt)\int_0^x\exp\lt(\int_0^z\frac{\wit b(u)}{\wit a(u)}\,du\rt)\frac{dz}{\wit a(z)}\, dx&<&+\iy\eqn
Taking into account the expressions for $\wit a$ and $\wit b$, 
we compute that
\bq
\fo x\in\RR_+,\qquad
\int_0^x\frac{\wit b(y)}{\wit a(y)}\,dy&=&\int_0^x\frac{a'(y) -b(y)+2a(y)k'(y)}{a(y)}\,dy\\
&=&\ln\lt(\frac{a(x}{a(0)}\rt)-c(x)+2(k(x)-k(0))\eq
so that
the l.h.s.\ of (\ref{II}) is proportional to 
the quantity
\bq
\int_0^{+\iy}\lt(\int_0^{x}
(\mu((-\iy,y]))^2\exp(-c(y))\, dy\rt) \frac1{
(\mu((-\iy,x]))^2}
\,\mu(dx)\eq
Since $0<\mu((-\iy,0])\leq \mu((-\iy,x])\leq 1$ for $x\in\RR_+$, the finiteness of the previous expression is equivalent to that of
$I_+$.
The only problem is that the coefficients $\wit a$ and $\wit b$ were required to be of  class $\cC^1$ by Ikeda and Watanabe.
But one can check directly in Section 3 of Chapter 6 of their book  \cite{MR1011252}
that the proof extends to the situation where the lack of regularity is restricted to 0,
where $\wit a$ is assumed to be continuous and positive and $\wit b$ locally bounded.
Alternatively, one can come back to the smooth situation in the following way.
Define
\bq
\fo x\in\RR\sqcup\{-\iy,+\iy\},\qquad \tau(V,x)&\df&\inf\{t\geq 0\st V_t=x\}\eq
As a consequence of the Markov property and of the symmetry of $V$, the fact that $\tau(V)$ is finite a.s., whatever the initial point, is equivalent to
\bqn{probas}
\lt\{
\begin{array}{rcl}
\PP[\tau(V,-2)\wedge \tau(V,2)<+\iy\vert V_0=0]&=&1\\
\PP[\tau(V,+\iy)<\tau(V,1)\vert V_0=2]&>&0
\end{array}\rt.
\eqn
The Girsanov transformation enables to see that the first of these conditions is true
as soon as $\wit a$ is continuous and positive on $[-2,2]$ and $\wit b$ is bounded on $[-2,2]$.
The second condition is not affected by modifications of $\wit a$ and $\wit b$ in $(-2,2)$.
So we can first apply Theorem 3.2 (3) of Chapter 6 of the book of Ikeda and Watanabe \cite{MR1011252} 
to  symmetric smoothings of $\wit a$ and $\wit b$ in $(-2,2)$ (this does not change the condition $I_+<+\iy$ either)
and next deduce the same conclusion for the original process $V$ via (\ref{probas}).\wwtbp
\par
Now we have at our disposal all the ingredients necessary to the proof of Proposition \ref{Itau+}.
So let us assume that $I_+<+\iy$.
\par
We begin by defining the following stopping times.
\bq
\wi\sigma_0&\df&\inf\{t\geq 0\st Y_t^*\geq 1\}\eq
By Corollary
\ref{convZ} we already know that $\lim_{t\ri+\iy} Y^*_t=+\iy$ so that $\wi\sigma_0$ is finite a.s.
Next consider
\bq
\wit\sigma_0&\df&\inf\{t> \wi\sigma_0\st Y_t^*=+\iy \hbox{ or } Y_t^*=0\}\eq
Since $I_+<+\iy$, we deduce from Lemma \ref{UY} and Proposition \ref{Uiy}  that
$\wit\sigma_0$ is finite a.s.\ and we have either $Y^*_{\wit\sigma_0}=+\iy$ or $Y^*_{\wit\sigma_0}=0$.
Indeed knowing the trajectory $Z^*_{[0,\wi\sigma_0]}$, the conditional probability that $Y^*_{\wit\sigma_0}=+\iy$
is bounded below by $\PP[\tau(U)<\tau(U,0)\vert U_0=Y^*_{\wi\sigma_0}]$,
where $\tau(U,0)=\inf\{t\geq 0\st U_t=0\}$.
Since the mapping $\RR_+\ni x\mapsto \PP[\tau(U)<\tau(U,0)\vert U_0=x]$ is non-decreasing, we get that 
\bq
\PP[Y^*_{\wit\sigma_0}=+\iy\vert \cF^*_{\wi\sigma_0}]\ \geq \ p_*\ \df\ \PP[\tau(U)<\tau(U,0)\vert U_0=1]\ >\ 0\eq
where $\cF^*_\sigma$ will designate the $\sigma$-field associated to the stopping time $\sigma$
in the filtration generated by the process $Z^*$: more explicitly $\cF^*_\sigma$ is generated by the piece of trajectory
$Z^*_{[0,\sigma]}$ (see e.g.\ Chapter 1 of the book \cite{MR1725357} of Revuz and Yor). It follows that $\PP[Y^*_{\wit\sigma_0}=+\iy]\geq p_*$.
If $Y^*_{\wit\sigma_0}=+\iy$, we set $N=0$ and otherwise the value of the random variable $N$
will be defined later on in the procedure.
Indeed, if $Y^*_{\wit\sigma_0}=0$, we consider
\bq
\wi\sigma_1&\df&\inf\{t> \wit\sigma_0\st Y_t^*= 1\}\\
\wit\sigma_1&\df&\inf\{t> \wi\sigma_1\st Y_t^*=+\iy \hbox{ or } Y_t^*=0\}\eq
These stopping times are again a.s.\ finite (still conditionally on $Y^*_{\wit\sigma_0}=0$).
If $Y^*_{\wit\sigma_1}=+\iy$, we set $N=1$. Note that as before,
\bq
\PP[Y^*_{\wit\sigma_1}=+\iy\vert \cF^*_{\wit\sigma_0},\, Y^*_{\wit\sigma_0}=0]& \geq & p_*\eq
The construction goes on similarly: if for some $n\in\NN$,
$\wit \sigma_n$ has been defined, we set $N=n$ if $Y^*_{\wit\sigma_n}=+\iy$
and the procedure stops. Otherwise,
namely if $Y^*_{\wit\sigma_n}=0$, we consider the a.s.\ finite random times
\bq
\wi\sigma_{n+1}&\df&\inf\{t> \wit\sigma_n\st Y_t^*= 1\}\\
\wit\sigma_{n+1}&\df&\inf\{t> \wi\sigma_{n+1}\st Y_t^*=+\iy \hbox{ or } Y_t^*=0\}
\eq
and we set $N=n+1$ if $Y^*_{\wit\sigma_{n+1}}=+\iy$.
The previous arguments show that
\bq
\PP[Y^*_{\wit\sigma_{n+1}}=+\iy\vert \cF^*_{\wit\sigma_n},\, Y^*_{\wit\sigma_n}=0]& \geq & p_*\eq
The validity of this property for all $n\in \ZZ_0^+$ implies that $N$ is stochastically bounded
below by a geometric random variable of parameter $1-p_*<1$:
\bq
\fo n\in  \ZZ_0^+,\qquad \PP[N\geq  n]&\leq & (1-p_*)^n\eq
In particular, $N$ is a.s.\ finite as well as $\tau^+=\wit \sigma_{N}$. This ends the proof of Proposition  \ref{Itau+}
and the finiteness assertion of Proposition \ref{Itau}. As explained in the introduction, this implies that $\tau$ is a strong stationary time for $X$,
once $X$ and $Z^*$ are intertwined through $\Lambda$.

\section{Intertwining}\label{SecI}

In the two previous sections, the process $Z^*$ has been studied in some details. It is time now to check that it can be intertwined
with the initial one-dimensional positive recurrent diffusion $X$.
\par\me
We begin by verifying  that the commutation relation (\ref{inter}) is satisfied with $L^*$ defined by
(\ref{Lstar}), (\ref{Lstar2}) and (\ref{Lstar3}). 
\begin{lem}\label{commrel}
For any $f\in\cC^2(\RR)$ such that $f$ and $L[f]$ belong to $\LL^1(\mu)$, we have
\bq
\fo z\in E^*\setminus (D^*\sqcup\{(-\iy,+\iy)\}),\qquad \Lambda[ L[f]](z)&=&L^* [\Lambda[f]](z)\eq
\end{lem}
\proof
A priori there are three situations to be considered
$z\in \rE$, $z\in\{-\iy\}\times\RR$ and $z\in \RR\times \{+\iy\}$.
We are to deal only with the first case, the other ones being similar (and even easier).
So let  $f\in\cC^2(\RR)$ be given (the integrability assumptions are needed only for 
$z\in\{-\iy\}\times\RR$ and $z\in \RR\times \{+\iy\}$ to insure the integrability of $f$ and $L[f]$ with respect to
$\mu$ on semi-infinite intervals).
For $z\df(x,y)\in\RR^2$ with $x<y$, we have
\bq
\Lambda[f](z)&=&\frac1{h(x,y)}\int_x^y f(u)\, \mu(du)\eq
where $h$ was defined in (\ref{h}). Taking into account Lemma \ref{evolset},
we get that
\bq
L^*[\Lambda[f]](z)&=&\frac{1}{h(z)}\wi L[F](z)\eq
where $\wi L$ was given in (\ref{wiL}) and where $F$ is the function defined on $E^*$ by 
\bqn{F}
\fo (x',y')\in E^*,\qquad F(x',y')&\df& \int_{x'}^{y'} f(u)\, \mu(du)\eqn
For $(x,y)\in \rE$, $\pa_x F(x,y) =-\mu(x)f(x)$ and $\pa_y F(x,y) =\mu(y)f(y)$, so that we get  that
\bqn{wiLF}
\nonumber\wi L[F](x,y)&=& (\sqrt{a(y)}\pa_y-\sqrt{a(x)}\pa_x)(\sqrt{a(y)}\mu(y)f(y)+\sqrt{a(x)}\mu(x)f(x))\\
\nonumber&&
-(a'(x)/2-b(x))\mu(x)f(x)+(a'(y)/2-b(y))\mu(y)f(y)\\
&=&
\sqrt{a(y)}\pa_y(\sqrt{a(y)}\mu(y)f(y))-\sqrt{a(x)}\pa_x(\sqrt{a(x)}\mu(x)f(x))\\
\nonumber&&
-(a'(x)/2-b(x))\mu(x)f(x)+(a'(y)/2-b(y))\mu(y)f(y)\\
\nonumber&=&
a(y)\mu(y)\pa_yf(y)-a(x)\mu(x)\pa_yf(x)- g(x)f(x)+g(y)f(y)
\eqn
where $g$ is the function defined by
\bq
\fo x\in \RR,\qquad g(x)&\df& \sqrt{a(x)}\pa_x(\sqrt{a(x)}\mu(x))+(a'(x)/2-b(x))\mu(x)\eq
Recalling the definition of $\mu$ given in (\ref{mu2}), we compute that $g$ is vanishes identically, so that
we obtain
\bq
\fo (x,y)\in\rE,\qquad 
L^*[\Lambda[f]](x,y)&=&\frac1{h(x,y)}\lt(a(y)\mu(y)\pa_yf(y)-a(x)\mu(x)\pa_yf(x)\rt)\eq
We turn now to the computation of  $\Lambda[L[f]]$ on $\rE$.
Note that 
$L$ can be factorized into
\bq
L\ \cdot\ = 
\ a\exp(-c)\pa(\exp(c)\pa\ \cdot)\ =\ \frac1{\mu}\pa(a\mu\pa\ \cdot)\eq
It follows that for all $f\in\cC^2(\RR)$ and $(x,y)\in\rE$, 
\bq
\int_x^y L[f](u)\,\mu(du)&=&\int_x^y \pa(a\mu\pa f)(u)\, du\\
&=&a(y)\mu(y) f'(y)-a(x)\mu(x) f'(x)\eq
The wanted commutation relation follows at once, on $\rE$.\wwtbp
\begin{rem}
If in the above proof $\wi L$ is replaced by
the generators $\check L$ or $\breve L$ defined respectively by (\ref{checkL}) and 
 (\ref{breveL}) (on $\rE$ and their natural extensions on $\RR\times\{+\iy\}$ and by $\{-\iy\}\times\RR$),
 the same computations are still valid. Indeed, remark that in (\ref{wiLF})
 the cross differentiation $\pa_x\pa_y$ vanishes, meaning that on $\rE$ and for the function $F$ defined
 in (\ref{F}),
 \bq
 \wi L [F]\ =\ \check L[F]\ =\ \breve L[F]\eq
 (simpler considerations are also valid on 
 $\{-\iy\}\times\RR\sqcup \RR\times \{+\iy\}$).
\\
The commutation relations $\Lambda L=\check L^*\Lambda$ and $\Lambda  L=\breve L^*\Lambda$
for the generators $\check L^*$ and $\breve L^*$ (defined respectively in (\ref{dagger}) and (\ref{ddagger})) are then also true,
because these operators are the $h$-transforms of  $\check L$ and $\breve L$, as it was mentioned in Remark \ref{dddagger}.
This justifies the assertions made after Proposition \ref{entrance} in the introduction.
\par
Even if some of the subsequent developments could
be extended to these generators, recall that their interest is limited, due to 
the observations made in Remark \ref{alpha}.
\end{rem}
\par
We are now going to lift the commutation relation of  Lemma \ref{commrel} to the level of the corresponding semi-groups. More precisely,
let $(P_t)_{t\geq 0}$ be the semi-group associated to $L$.
From a probabilistic point of view, it is constructed in the following way.
For any $x\in \RR$, consider $(X_t)_{t\in\RR}$
the solution starting from $x$ of the s.d.e.\ 
\bqn{X}
dX_t &=& b(X_t)dt +\sqrt{2a(X_t)}dB_t\eqn
where $(B_t)_{t\geq 0}$ is a standard Brownian motion.  
Then for any $t\geq 0$ and any bounded and continuous mapping $f$ on $\RR$, we have
\bq
P_t[f](x)&=&\EE_x[f(X_t)
]\eq
\par\sm
The semi-group $(P^*_t)_{t\geq 0}$ can be constructed similarly. For $z\in E^*\cap\RR\times\RR$, consider the process $Z^*$
starting from $z$ defined in Proposition \ref{entrance} or \ref{facile}, depending if $z\in D^*$ or not
(if $z=(-\iy,+\iy)$, $Z^*$ stays forever at $(-\iy,+\iy)$). For 
$z\in\{-\iy\}\times\RR$ or $z\in \RR\times \{+\iy\}$, $Z^*$ is constructed as explained in the proof of Proposition \ref{facile}.
Then for any
 $t\geq 0$ and any bounded and continuous
 mapping $f$ on $E^*$, we take
\bq
P_t^*[f](z)&=&\EE_z[f(Z^*_t)
]\eq
\par
\begin{pro}\label{LPPL}
Assume that $X$ is positive recurrent.
Then for all $T\geq 0$ and all bounded and continuous  function $f$ on $\RR$, we have
\bq
\fo z\in E^*,
\qquad \Lambda[ P_T[f]](z)&=&P_T^* [\Lambda[f]](z)\eq
\end{pro}
Formally, writing $P_t=\exp(tL)$ and $P_t^*=\exp(tL^*)$, the deduction of these commutation relations from their infinitesimal version given in Lemma \ref{commrel} may seem clear. Nevertheless a direct rigorous justification does not seem so obvious (see Remark \ref{Hilbert} below). 
We found it preferable to follow a recurrent idea in the study of semi-groups à la
 Bakry \cite{MR1307413} and Ledoux \cite{MR1600888}.
\par
\proof
It consists in investigating the evolution of
\bq
[0,T]\ni t&\mapsto& P^*_{t}[\Lambda[P_{T-t}[f]]]\eq
for given $T>0$ and first for $f\in\cCc^\iy(\RR)$.\par
We begin by recalling how to exploit the martingale property of  $Z^*$.
A function defined on $E^*$ is said to be $\cC^2$ if it is continuous on $E^*$ and if it is
$\cC^2$ on $\rE$, on $\{-\iy\}\sqcup\RR$ and on $\RR\sqcup\{+\iy\}$.
Similarly, a continuous function defined on $\RR_+\times E^*$ is said to be $\cC^{1,2}$ if it is $\cC^1$ with respect with the first variable in $\RR_+$
and $\cC^2$ with respect to the second variable in $E^*$, the corresponding partial derivatives
being continuous on $\RR_+\times\rE$, on $\RR_+\times(\{-\iy\}\times\RR)$ and on $\RR_+\times(\RR\times\{+\iy\})$.
Denote by $\cCb^{1,2}(\RR_+\times E^*)$ the set of such functions $F$ which are furthermore bounded, as well as the mapping $\RR_+\times(E^*\setminus(D^*\sqcup \{(-\iy,+\iy)\}))\ni(t,z)\mapsto\pa_tF(t,z)+L^*[F(t,\cdot)](z)$.
Let us prove that for any $z\in E^*\setminus D^*$, $t\geq 0$ and $F\in \cCb^{1,2}(\RR_+\times E^*)$,
\bqn{Pstart}
\EE_z[F(t\wedge \tau^*,Z^*_{t\wedge \tau^*})
]&=& F(0,z)+\EE_z\lt[\int_0^{t\wedge \tau^*} \pa_sF(s,Z_s^*)+L^*[F(s,\cdot)](Z^*_s)\,ds\rt]\eqn
First we treat the case  where $z=(x,y)\in\rE$ and we replace $\tau^*$ by $\tau^\dagger$ which was defined in \eqref{taudagger}.
Indeed, for $n\in\NN$ large enough, say $n\geq n_0$, where $n_0\in\NN$ is  such that $y-x>1/n_0$, consider
\bq
\tau_1(n)&\df&\inf\{t\geq 0\st X^*_t< -n\}\\
\tau_2(n)&\df&\inf\{t\geq 0\st Y^*_t > n\}\\
\tau_3(n)&\df&\inf\{t\geq  0\st Y^*_t-X^*_t<1/n\}\\
\tau^{\dagger}(n)&\df&\min( \tau_1(n),\tau_2(n),\tau_3(n))
\eq
where $Z^*=(X^*_t,Y_t^*)_{t\geq 0}$.
The sequence $(\tau^{\dagger}(n))_{n\geq n_0}$ is a localizing sequence for $Z^*$ on the random time interval $[0,\tau^{\dagger})$,
in the sense that
\bq
\tau^\dagger&=&\lim_{n\ri\iy}\tau^{\dagger}(n)\eq
and for any $F\in \cCb^{1,2}(\RR_+\times E^*)$, we can write
\bq
\fo t\geq 0,\qquad
F(t\wedge \tau^\dagger, Z^*_{t\wedge \tau^\dagger})&=&F(0,z)+
\int_0^{t\wedge \tau^\dagger} \pa_sF(s,Z_s^*)+L^*[F(s,\cdot)](Z^*_s)\,ds+M_{t}\eq
where for any $n\geq n_0$,
the process $(M_{t\wedge \tau^{\dagger}(n)})_{t\geq 0}$ is a martingale starting from 0.\\
So taking expectations, we end up with
\bq
\EE_z[F(t\wedge \tau^\dagger(n), Z^*_{t\wedge \tau^\dagger(n)})]&=&F(0,z)+
\EE_z\lt[\int_0^{t\wedge \tau^\dagger(n)} \pa_sF(s,Z_s^*)+L^*[F(s,\cdot)](Z^*_s)\,ds\rt]\eq
Our boundedness and continuity assumptions on $F$ enable to use the bounded convergence theorem
to get 
\bqn{Pstartau}
\EE_z[F(t\wedge \tau^\dagger, Z^*_{t\wedge \tau^\dagger})]&=&F(0,z)+
\EE_z\lt[\int_0^{t\wedge \tau^\dagger} \pa_sF(s,Z_s^*)+L^*[F(s,\cdot)](Z^*_s)\,ds\rt]\eqn\par
Recall from Lemma \ref{mart} that if $\tau^\dagger<+\iy$, then $Z^*_\dagger$ belongs
to $\{(-\iy,+\iy)\}\sqcup( \{-\iy\}\times \RR)\sqcup(\RR\times\{+\iy\})$.
Note also that if $\tau^\dagger<+\iy$ and $Z^*_{\tau^\dagger}=(-\iy,+\iy)$, then $\tau^*=\tau^\dagger$, so 
\bq\EE_z[F(t\wedge \tau^\dagger, Z^*_{t\wedge \tau^\dagger})\un_{\{\tau^\dagger\leq t,\, Z^*_{\tau^\dagger}= (-\iy,+\iy)\}}]&=&
\EE_z[F(t\wedge \tau^*, Z^*_{t\wedge \tau^*})\un_{\{\tau\dagger\leq t,\, Z^*_{\tau^\dagger}= (-\iy,+\iy)\}}]
\eq
Thus to prove \eqref{Pstart}, taking into account the
strong Markov property at time $t\wedge\tau^\dagger$ (which is true by construction of $Z^*$),
it is sufficient to see that for all $z\in ( \{-\iy\}\times \RR)\sqcup(\RR\times\{+\iy\})$,
\bq
\EE_z[F(t\wedge \tau^*, Z^*_{t\wedge \tau^*})]&=&F(0,z)+
\EE_z\lt[\int_0^{t\wedge \tau^*} \pa_sF(s,Z_s^*)+L^*[F(s,\cdot)](Z^*_s)\,ds\rt]\eq
This is immediate, following a localization procedure similar to that leading to \eqref{Pstartau}.
\\
Since $(-\iy,+\iy)$ is absorbing, if $\tau^*<t$ we can write for $F\in\cCb^{1,2}(\RR_+\times E^*)$,
\bq
F(t,Z_t^*)&=& F(t,(-\iy,+\iy))\\
&=&F(\tau^*,(-\iy,+\iy))+\int_{\tau^*}^t \pa_sF(s,(-\iy,+\iy))\,ds\\
&=&F(\tau^*,Z^*_{\tau^*})+\int_{\tau^*}^t \pa_sF(s,Z^*_{s})\,ds\eq
so that, recalling the Dirichlet condition for $L^*$ at $(-\iy,+\iy)$, \eqref{Pstartau} can be transformed into
\bq
\EE_z[F(t,Z^*_{t})]&=&
 F(0,z)+\EE_z\lt[\int_0^{t} \pa_sF(s,Z_s^*)+L^*[F(s,\cdot)](Z^*_s)\,ds\rt]\eq
 namely in semi-group notations,
 \bqn{Pstart2}
 P^*_t[F(t,\cdot)](z)= F(0,z)+\int_0^tP_s^*[\pa_sF(s,\cdot)+L^*[F(s,\cdot)]](z)\, ds\eqn
\par
Let $T>0$ and $f\in\cCc^\iy(\RR)$ be fixed, we want to apply the previous considerations with the function $F$ defined
on $[0,T]\times E^*$ by
\bq
\fo (t,z)\in [0,T]\times E^*,\qquad F(t,z)&\df& \Lambda[P_{T-t}[f]](z)\eq\par
Since $\RR_+\times\RR\ni(t,x)\mapsto P_t[f](x)$ is well-known to be smooth, it is clear that $F$ is $\cC^{1,2}$.
Furthermore, recall that the semi-group $(P_t)_{t\geq 0}$ can be extended into a self-adjoint continuous semi-group on $\LL^2(\mu)$,
whose generator is the Friedrich extension of $L$ on $\LL^2(\mu)$.
It follows that the relation $\pa_t P_t[f]=LP_t[f]$ is satisfied in the usual sense and in $\LL^2(\mu)$ and we get
\bqn{patF}
\fo t\in[0,T],\,\fo z\in E^*,\qquad \pa_tF(t,z)&=& -\Lambda[L[P_{T-t}[f]]](z)\eqn
Since  the mapping $\RR\ni x\mapsto P_{T-t}[f](x)$ is $\cC^2$ and
\bq \mu[\lve P_{T-t}[f]\rve]&\leq & \mu[\lve f\rve]\\
\mu[\lve L[P_{T-t}[f]]\rve]
&=&\mu[\lve P_{T-t}[L[f]]\rve]
\ \leq\  \mu[\lve L[f]\rve]\eq
 we are in position to apply
Lemma \ref{commrel} (with $f$ replaced by $P_{T-t}[f]$) to get that in the r.h.s.\ of \eqref{patF},
we can replace 
$\Lambda[L[P_{T-t}[f]]](z)$ by $L^*[\Lambda[P_{T-t}[f]]](z)$,
at least for $z\in E^*\setminus (D^*\sqcup\{(-\iy,+\iy)\})$.
Thus we get that
\bq
\fo t\in[0,T],\,\fo z\in E^*\setminus (D^*\sqcup\{(-\iy,+\iy)\}),\qquad\pa_tF(t,z)+L^*[F(t,\cdot)](z)&=&0\eq
This relation is also true for $z=(-\iy,+\iy)$. Indeed, due to the fact that $X$ is positive recurrent, we get
\bq
\fo t\geq 0,\qquad F(t,(-\iy,+\iy))\ =\ \mu[P_t[f]]\ =\ \mu[f]\eq
so that
\bqn{posrec}
\pa_tF(t,(-\iy,+\iy))&=&0\eqn
In particular it is licit to apply \eqref{Pstart2} (for $t\in[0,T]$) to get
\bq
\fo z\in E^*,\qquad P^*_T[F(T,\cdot)](z)&=& F(0,z)\eq
which is just the conclusion stated in the proposition, at least for $f\in\cCc^{\iy}(\RR)$.
To extend it to any bounded and continuous function $f$, note that for any fixed $T\geq 0$ and $z\in E^*$,
the mappings
\bq
\cR\ni A\mapsto \Lambda[P_T[\un_A]](z)&\hbox{ and }& \cR\ni A\mapsto P_T^*[\Lambda[\un_A]](z)\eq
($\cR$ stands for the $\sigma$-algebra of Borelian subsets of $\RR$) define  two probability measures. Because they
coincide on every $f\in\cCc^{\iy}(\RR)$, they must  be equal.\wwtbp
\begin{rem}
The assumption that $X$ is positive recurrent is really necessary for the previous result.
Indeed, there exists generators $L$ satisfying \eqref{mfini}
but  not \eqref{JJ}. For the associated semi-group $(P_t)_{t\geq 0}$, for any time $T>0$ and any point $x\in \RR$ we have
$P_T[\un](x)<1$. As a consequence, for any $T>0$ and $z\in E^*$, $\Lambda[P_T[\un]](z)<1$, while by construction
$P_T^*[\Lambda[\un]](z)=P_T^*[\un](z)=1$.\\
In the above proof, the positive recurrence of $X$ is encapsulated in \eqref{posrec}.
\end{rem}
\begin{rem}\label{Hilbert}
When $F$ doesn't depend on the time variable, \eqref{Pstart2} writes down under the familiar form
\bq
\pa_t P_t^*[F]&=&L^*[P_t^*[F]]\eq
But from an analytical point of view, it is not clear a priori in which Banach space one should interpret this evolution equation
 to deduce the semi-group $(P^*_t)_{t\geq 0}$ from $L^*$. 
 If we were to work with the elliptic generator $\breve L^*$ defined in \eqref{ddagger}, there is a natural $\LL^2$ Hilbert setting.
 Indeed, let $\wi \eta$ be the $\sigma$-finite measure on $\RR$ whose density with respect to the Lebesgue measure is 
 $\exp(-c)$. The generator $\breve L$ given in \eqref{breveL} is then symmetric with respect to the measure $\eta$
 which coincides on
 $E^*$ with the restriction of $(\wi \eta+\delta_{-\iy}+\delta_{+\iy})^{\otimes 2}$.
 Since  $\breve L^*$ corresponds to the $h$-transform of $\breve L$, it is symmetric relatively to the measure $\nu$
 admitting $h^2$ as density with respect to $\eta$.
 Thus the relations $\breve P^*_t=\exp(t\breve L^*)$, for $t\geq 0$, could be given a meaning in $\LL^2(\nu)$.
 Heuristically, the intertwining between $L$ and $L^*$ can be seen as ``weak conjugation relation'' between them,
 so we can
 expect that $L^*$ is equally reversible with respect to some $\sigma$-finite measure on $E^*$. Unfortunately we have not been
 able to find it and in addition
 we have no idea about possible quasi-invariant measures of $L^*$. Nevertheless, we believe that this subject really deserves to be investigated further, especially from a quantitative point of view.  An initiation of this program in a very particular case is presented in the next section.
\end{rem}
\par
Proposition \ref{LPPL} is the main technical point to get the intertwined coupling  of $X$ with $Z^*$. Indeed, we can follow the construction
of Diaconis and Fill \cite{MR1071805} by applying it to skeletons of $X$ with $Z^*$. Passing to the limit in the latter approximations will enable us to justify the arguments given before the statement of Proposition \ref{Itau} in the introduction.
\par\sm
Let be given $m_0$ and $m_0^*$ two probability measures respectively on $\RR$ and $E^*$ such that $m_0^*\Lambda=m_0$.
We want to construct an intertwining of $X$ with $Z^*$ whose initial distribution is described by
$\eta_0(dx,dz^*)\df m_0^*(dz^*)\Lambda(z^*,dx)$ (in particular the laws of $X_0$ and $Z_0^*$ are respectively $m_0$ and $m_0^*$).
For fixed $N\in\NN$, define a discrete time Markov chain $(\bar X^{(N)}_{n2^{-N}}, \bar Z^{(N,*)}_{n2^{-N}})_{n\in\ZZ_+}$,  intertwined through $\Lambda$,
in the following way: its initial distribution is $\eta_0$ and its transition kernel $Q^{(N)}$ is given by
\bq
Q^{(N)}((x,z^*), d(\wi  x,\wi  z^*))&\df&P_{2^{-N}}(x,d\wi  x)P^*_{2^{-N}}(z^*,d\wi  z^*)\frac{\Lambda(\wi  z^*,d\wi  x)}{\tr_{2^{-N}} (z^*,d\wi  x)}\eq
(from $(x,z^*)\in\RR_+\times E^*$ to the infinitesimal neighborhood $d(\wi  x,\wi  z^*)$ of  $(\wi x,\wi z^*)\in\RR_+\times E^*$), where the last ratio is the Radon-Nikodym derivative of the measure $\Lambda(\wi  z^*,d\wi  x)$
with respect to $\tr_{2^{-N}} (z^*,d\wi  x)\df (P^*_{2^{-N}}\Lambda)(z^*,d\wi  x)= (\Lambda P_{2^{-N}})(z^*,d\wi  x)$.
One would have remarked that due to Propositions \ref{entrance} and \ref{facile}, for any $z^*\in E^*$, $\tr_{2^{-N}} (z^*,\,\cdot\,)$
is equivalent to the Lebesgue measure. So except if $\wi  z^*$ corresponds to a singleton,
we have $\Lambda(\wi  z^*,\,\cdot\,)\ll \tr_{2^{-N}} (z^*,\,\cdot\,)$. But $P^*_{2^{-N}}(z^*,d\wi  z^*)$-a.s.\ 
$\wi  z^*$ does not correspond to a singleton, so 
$Q_{2^{-N}}$ is indeed a transition kernel (not only a sub-Markovian kernel).
\\
The computations of Diaconis and Fill \cite{MR1071805} can then be adapted to this setting, because of the structure of the initial distribution and of Proposition \ref{LPPL},
to show that the Markov chain $(\bar X^{(N)}_{n2^{-N}}, \bar Z^{(N,*)}_{n2^{-N}})_{n\in\ZZ_+}$ thus constructed satisfies the following properties: 
\bqn{intert1}
\hbox{
$( \bar X^{(N)}_{n 2^{-N}})_{n\in\ZZ_+}$ and $(X^{(N)}_{n 2^{-N}})_{n\in\ZZ_+}$ have the same law}\eqn
\bqn{intert2}
\hbox{
$( \bar Z^{(N,*)}_{n 2^{-N}})_{n\in\ZZ_+}$ and $(Z^{(N,*)}_{n 2^{-N}})_{n\in\ZZ_+}$ have the same law}\eqn
\bqn{intert3}
\hbox{$\fo m\in\ZZ_+$, the conditional law of $\bar X^{(N)}_{m2^{-N}}$
knowing $\bar Z^{(N,*)}_{0},\, \bar Z^{(N,*)}_{2^{-N}}, \, ..., \bar Z^{(N,*)}_{m2^{-N}}$ is $\Lambda(\bar Z^{(N,*)}_{m2^{-N}},\,\cdot\,)$}\eqn
\bqn{intert4}
\nonumber&\hbox{$\fo m\in\ZZ_+$, the conditional law of $(\bar Z^{(N,*)}_{0},\, \bar Z^{(N,*)}_{2^{-N}}, \, ..., \bar Z^{(N,*)}_{m2^{-N}})$
knowing $( \bar X^{(N)}_{n 2^{-N}})_{n\in\ZZ_+}$}&\\
&\hbox{depends only on  $\bar X^{(N)}_{0},\, \bar X^{(N)}_{2^{-N}}, \, ..., \bar X^{(N)}_{m2^{-N}}$}&\eqn
\par
Next we embed the 
Markov chain $(\bar X^{(N)}_{n2^{-N}}, \bar Z^{(N,*)}_{n2^{-N}})_{n\in\ZZ_+}$ into the (time-inhomogeneous) Markov process
$(\bar X^{(N)}, \bar Z^{(N,*)})\df(\bar X^{(N)}_{t}, \bar Z^{(N,*)}_{t})_{t\in\RR_+}$,
by taking
\bq
\fo t\geq 0,\qquad (\bar X^{(N)}_{t}, \bar Z^{(N,*)}_{t})\df (\bar X^{(N)}_{\lfloor t2^{N}\rfloor2^{-N}}, \bar Z^{(N,*)}_{\lfloor t2^{N}\rfloor2^{-N}})\eq
where $\lfloor\cdot\rfloor$ stands for the integer part.
\begin{pro}\label{XZstar}
The sequence of the laws of $(\bar X^{(N)}, \bar Z^{(N,*)})$, for $N\in\NN$, on the Skorokhod space $\DD(\RR_+,\RR\times E^*)$,
is the relatively compact. We can thus extract a subsequence converging to a probability measure $\PP$ which is necessarily  supported by the 
set of continuous trajectories. 
The canonical coordinate process $(\bar X_{t}, \bar Z^{*}_{t})_{t\in\RR_+}$ is a coupling of $X$ with $Z^*$ satisfying for all $ t\in\RR_+$, 
\bqn{intert5}
\hbox{the conditional law of $\bar X_{t}$
knowing $\bar Z^{*}_{[0,t]}$ is $\Lambda(\bar Z^{*}_{t},\,\cdot\,)$}\eqn
\bqn{intert6}
&\hbox{the conditional law of $\bar Z^{*}_{[0,t]}$
knowing $ \bar X$ depends only on  $\bar X_{[0,t]}$}&\eqn
\end{pro}
\proof
Using traditional properties of the Skorokhod topology on the Polish space $\DD(\RR_+,\RR\times E^*)$ (see for instance the book \cite{MR0233396} of Billingsley),
we deduce from \eqref{intert1} and \eqref{intert2} that the laws of $\bar X^{(N)}$ and $\bar Z^{(N,*)}$ converge
respectively toward those of $X$ and $Z^*$.
This observation implies without difficulty the first three assertions of the above proposition.
For the last two ones, note that as consequences of \eqref{intert3} and \eqref{intert4}, we have
\bq
\hbox{$\fo t\geq 0$, the conditional law of $\bar X^{(N)}_{t}$
knowing $\bar Z^{(N,*)}_{[0,t]}$ is $\Lambda(\bar Z^{(N,*)}_{t},\,\cdot\,)$}\eq
\bq
\nonumber&\hbox{$\fo t\geq 0$, the conditional law of $\bar Z^{(N,*)}_{[0,t]}$
knowing $\bar X^{(N)}$
depends only on  $\bar X^{(N)}_{[0,t]}$}&\eq
The deduction of  \eqref{intert5} and \eqref{intert6} is then a standard exercise on conditional expectations:
use on one hand that the $\sigma$-algebra generated by $\xi_{[0,t]}$, where $t\geq 0$ and $\xi$ is either $\bar X$ or $\bar Z^*$,
is the same as that generated by  mappings of the form $F(\xi_{t_1}, ..., \xi_{t_r})$,
where $r\in\NN$, $t_1,..., t_r$ are dyadic numbers satisfying $0\leq t_1<\cdots <t_r\leq t$ and $F$ is a bounded and continuous function on either
$\RR^r$ or $(E^*)^r$, and on the other hand that such mappings are $\PP$-a.s.\ continuous.
\wwtbp
\begin{rem}
Pal and Shkolnikov \cite{2013arXiv1306.0857P} investigated the existence
of intertwinings between diffusion semi-groups whose generators are appropriately linked by a Markov kernel.
Unfortunately the assumptions of their Theorem 3 do not cover our situation, essentially due the lack of ellipticity of $Z^*$.
The intertwining of $L$ with $ \breve L^*$ (defined in
\eqref{ddagger}) is  more amenable to their conditions, after 
looking at $X$ through the chart
$\RR\ni s\mapsto \int^s_0a^{-1/2}(u)\,du$ (and correspondingly for $Z^*$).
Nevertheless it would still remain to check their boundary conditions.
In the approach presented above, we escaped the corresponding  delicate description of what happens to the intertwined process
$(X,Z^*)$ when $X$ enters in contact with one of the boundaries of $Z^*$ by resorting to
the computations of Diaconis and Fill \cite{MR1071805} applied to the skeleton chains.
\end{rem}
\begin{rem}
The previous intertwinings of the skeleton chains are in general not compatible: 
it is not true that for all $N\in\NN$, $(\bar X^{(N+1)}_{n2^{-N}},\bar Z^{(N+1,*)}_{n2^{-N}})_{n\in\NN} $ has the same law as $(\bar X^{(N)}_{n2^{-N}},\bar Z^{(N,*)}_{n2^{-N}})_{n\in\NN}  $. 
\end{rem}
\par
Proposition \ref{XZstar} enables to prove the direct part of Theorem \ref{t1} along the arguments given before
Proposition \ref{entrance}. To end this section, we show the converse implication, by considering the diffusion $X$ 
whose initial distribution is $\mu$ conditioned to be on $\RR_-$ (namely  $\Lambda((-\iy,0),\cdot)$, the cases where
the initial distribution is $\Lambda((-\iy,x),\cdot)$ or $\Lambda((x,+\iy),\cdot)$, for some $x\in\RR$, can be treated similarly).
\par
In this situation the process $Z^*$ has the form $(-\iy, Y^*)$, where $Y^*$ is the solution starting from 0 of the s.d.e.\
\bq
dY^*_t&=&\lt(a'(Y^*_t)-b(Y^*_t)+2\frac{\sqrt{a(Y^*_t)}\mu(Y^*_t)}{\mu((-\iy,Y^*_t])}\sqrt{a(Y^*_t)}\rt)dt +\sqrt{2a(Y^*_t)}\,dB_t
\eq
From Corollary \ref{convZ}, we know a priori that $\lim_{t\ri+\iy}Y_t^*=+\iy$.
From Section \ref{Secei}, the boundary  $+\iy$ will be reached in finite time (a.s.) if and only if $I_+<+\iy$.
The reaching time of $+\iy$ by $Y^*$ is indeed
 the random time 
$\tau^*$ defined in (\ref{taustar}).
If we assume that $X$ admits a strong stationary time and if we show that such a strong stationary time is stochastically
larger than $\tau^*$, we would then get that $I_+<+\iy$. 
Symmetrically we would prove that the existence of a strong stationary time for $X$ starting from $\Lambda((0,+\iy),\cdot)$
implies that $I_-<+\iy$ and the converse part of Theorem \ref{t1} will be shown.
Thus according to (\ref{sharp}), it remains to check that
\begin{lem}
Under the previous assumption on $X$,  we have
\bq\fo t\geq 0,\qquad
 \mathfrak{s}(\cL(X_t),\mu)&= & \PP[\tau^* >  t]\eq
\end{lem}
The following arguments are an adaptation to our setting of Remark 2.39 of Diaconis and Fill \cite{MR1071805}.
\proof
Consider the intertwining of $X$ and $Z^*=(-\iy,Y^*)$ obtained in 
Proposition \ref{XZstar}. It follows that for all $t\geq 0$,
\bq
\cL(X_t)&=&\EE[\Lambda((-\iy,Y^*_t),\,\cdot\,)]\eq
In particular, we get that
\bq
 \mathfrak{s}(\cL(X_t),\mu)&= & \sup_{x\in\RR}
 \EE\lt[1-\frac{ d\Lambda((-\iy,Y^*_t),\,\cdot\,)}{d\mu}(x)\rt]\\
 &=&1-\inf_{x\in\RR}\EE\lt[\frac{ d\Lambda((-\iy,Y^*_t),\,\cdot\,)}{d\mu}(x)\rt]\eq
 The above Radon-Nikodym derivative is easy to compute: for all $x\in\RR$,
 \bq
 \frac{ d\Lambda((-\iy,Y^*_t),\,\cdot\,)}{d\mu}(x)&=&\frac{1}{\mu((-\iy,Y^*_t))}\un_{(-\iy,Y_t^*)}(x)\eq
 Note that the r.h.s.\ is non-increasing as a function of $x\in\RR$, so the same is true
 of the expression $\EE\lt[\frac{ d\Lambda((-\iy,Y^*_t),\,\cdot\,)}{d\mu}(x)\rt]$
 and we get
 \bq
  \mathfrak{s}(\cL(X_t),\mu)&= &1-\lim_{x\ri+\iy}\EE\lt[\frac{ d\Lambda((-\iy,Y^*_t),\,\cdot\,)}{d\mu}(x)\rt]\\
  &=&1-\PP[Y^*_t=+\iy]\\
  &=&\PP[Y^*_t<+\iy]\\
  &=&\PP[\tau^*<t]\eq
  \wwtbp

\section{On the Ornstein-Uhlenbeck counter-example}

In the study of convergence to equilibrium for diffusions, the  Ornstein-Uhlenbeck process is a benchmark,
in particular due to its Gaussian feature which enables explicit computations. 
Unfortunately, it is in some sense at the ``outside boundary" of the domain of application of the approach presented before.
We will see here how the method can nevertheless  be adapted to recover sharp informations.
\par\me
The Ornstein-Uhlenbeck process corresponds to the choice  in \eqref{L} of $a\equiv 1$ and $b(x)=-x$, for all $x\in\RR$.
The associated reversible measure is the centered and standard Gaussian distribution $\gamma$
whose density is given by $\gamma(x)=\exp(-x^2/2)/\sqrt{2\pi}$, for all $x\in\RR$.
A traditional integration by part leads to 
\bq
\gamma([x,+\iy))&\sim& \frac{\gamma(x)}{x}\eq
as $x$ goes to $+\iy$, 
so we get that the second integral of the l.h.s.\ of \eqref{IU} is infinite.
Theorem \ref{t1} then asserts that there exists initial distributions for which it is not possible to construct strong stationary times for  the associated
process $X$. 
Indeed, this is true as soon as the initial distribution $m_0$ has a compact support.
To see it, let us recall how the law $\cL(X_t)$ is easily computed in this situation: since $X$ 
satisfies the s.d.e.
\bq
\fo t\geq 0,\qquad dX_t&=&-X_t\,dt+\sqrt{2} dB_t\eq
(where $B=(B_t)_{t\geq 0}$ is a standard  Brownian motion), the variation of parameters method gives us:
\bq
X_t&=&\exp(-t)X_0+\sqrt{2}\int_0^t\exp(s-t)\, dB_s\eq
It follows that $m_t\df\cL(X_t)$ is the convolution of $\cL(\exp(-t)X_0)$ with $\gamma_{1-\exp(-2t)}$, the centered Gaussian distribution
of variance $1-\exp(-2t)=2\int_0^t\exp(2(s-t))\, ds$. Thus if $m_0$ has compact support, we get that for any fixed $t>0$, the separation discrepancy of
$m_t$ with $\gamma$ is one:
\bq
\mathfrak{s}(m_t,\gamma)\ =\ \lim_{\lve x\rve \ri +\iy}1-\frac{dm_t}{d\gamma}(x)\ =\ 1\eq
(a similar reasoning, considering only the limit at $-\iy$ or $+\iy$, would lead to the same conclusion
if the support of $m_0$ is bounded below or above: this enables to include the initial distributions 
considered for the reverse part of Theorem \ref{t1}).
The bound \eqref{sharp} then implies that there is no strong stationary time for $X$.\par
To simplify the presentation, we will assume that the initial distribution is the Dirac mass at 0.
We deduce from the above considerations that for any $t>0$, $\cL(X_t)=\gamma_{1-\exp(-2t)}$.
In particular $\cL(X_t)$ converges toward $\gamma$ in total variation. Let us check that the exponential rate for this
convergence is 2:
\begin{lem}\label{taux2}
We have
\bq
\lim_{t\ri+\iy}\frac1{t}\ln(\lVe m_t-\gamma\rVe_{\mathrm{tv}})&=&-2\eq
\end{lem}
\proof
By one of the characterization of the total variation norm, we have for all $t>0$,
\bqn{tv}
\lVe m_t-\gamma\rVe_{\mathrm{tv}}&=&\int (f_t-1)_+\,d\gamma\eqn
where $f_t$ is the Radon-Nikodym derivative of $m_t$ with respect to $\gamma$.
We compute 
that
\bq
\fo x\in\RR,\qquad
f_t(x)&=&(1-e^{-2t})^{-1/2}\exp\lt(-\frac{e^{-2t}x^2}{2(1-e^{-2t})}\rt)\eq
and 
we deduce that
\bq
f_t(x)\geq 1&\Longleftrightarrow& \lve x\rve\leq  x_t
\ \df\ \sqrt{(e^{2t}-1)\ln(1-e^{-2t})}
\eq
The quantity $x_t$ converges toward 1 when $t$ goes to infinity.
A simple expansion of the expression $f_t(x)-1$ then leads to
\bq
\lVe m_t-\gamma\rVe_{\mathrm{tv}}&=&2\int_0^{x_t} f_t(x)-1\, \gamma(dx)\\
&\sim& e^{-2t}\int_0^11-x^2\, \gamma(dx)\eq
for large $t>0$. The announced result follows at once.\wwtbp
\begin{rem}
The logarithmic Sobolev constant associated to $L$ is 4, so starting from any initial
distribution $m_0$ such that 
the relative entropy of $m_t$ with respect to $\gamma$ is finite for some $t\geq 0$,
we get that the exponential rate of converge in the relative entropy sense is at least 4.
Using next Pinsker's inequality, we recover that the above exponential rate of convergence
in total variation is at least 2.
For this traditional approach, see for instance the book \cite{MR2002g:46132} of 
{An{\'e},  Blach{\`e}re,               Chafa{\"{\i}},  Foug{\`e}res,  Gentil,
            Malrieu, Roberto and Scheffer}.
The Ornstein-Uhlenbeck process is also critical for the use of the logarithmic Sobolev inequalities method,
but it is in the ``interior boundary" of the domain of application.
              \end{rem}
\par
Let us show how to recover this exponential rate 2 for  the convergence in total variation by using strong (non-stationary) times.
So the emphasis is in testing the method, not in the result itself. It will also enable us to illustrate on this example
the directions suggested by Remark \ref{Hilbert}.
\par\sm
We begin by noting that the construction of the process $Z^*=(X^*,Y^*)$ made in Section \ref{dotdp}
is still valid. By symmetry and since we are considering $Z^*_0=(0,0)$,
we have that $X^*=-Y^*$. It comes from the fact that $Z^*$ and $(-Y^*,Y^*)$ satisfy the same well-posed martingale problem.
The diffusion $Y^*$ is given as the solution starting from 0 (which is an entrance boundary for $Y^*$)
of the s.d.e.
\bqn{Ystar}
\fo t> 0,\qquad
dY^*_t&=&\lt(Y^*_t+g(Y_t^*)\rt)dt +\sqrt{2}\,dB_t
\eqn
where as usual $B\df(B_t)_{t\geq 0}$ is a standard  Brownian motion, and where $g$ is the mapping
defined by
\bqn{g}
\fo y>0,\qquad g(y)&\df& 2\frac{\gamma(y)}{\gamma([0,y])}\eqn
The coupling of $X$ and $Y^*$ constructed in Section \ref{SecI} is equally valid. We deduce that any stopping time for $Y^*$
is a strong time for $X$. For any $M>0$, we are particularly interested in the following stopping time
\bq
\tau^*_M&\df& \inf\{t\geq 0\st Y^*_t=M\}\eq
It has the property that $\tau^*_M$ and $X_{\tau^*_M}$ are independent and that $X_{\tau^*_M}$
is distributed according to $\gamma_{[-M,M]}$, the conditioning of $\gamma$ on the interval $[-M,M]$.
The interest of the independence of the time and the position appears in the proof of
\begin{lem}\label{tvtv}
For all $t\geq 0$ and $M>0$, we have
\bq
\lVe m_t-\gamma\rVe_{\mathrm{tv}}&\leq & \PP[\tau^*_M>t]+\lVe\gamma_{[-M,M]}-\gamma\rVe_{\mathrm{tv}}\eq
\end{lem}
\proof
An equivalent formulation to \eqref{tv} of the total variation
is given by
\bqn{tv2}
\lVe m_t-\gamma\rVe_{\mathrm{tv}}&=&\frac12\sup_{\lVe f\rVe_{\iy}=1}\EE[f(X_t)]-\gamma[f]\eqn
where the supremum is taken over all measurable functions $f$ taking values in $[-1,1]$.
\\
Let $\cF_{\tau^*_M}$ be the $\sigma$-field generated by the piece of trajectory of the intertwined process
$(X,Y^*)$ up to time $\tau_M^*$. It is in fact generated by $X_{[0,\tau^*_M]}$ and some randomness independent
from the whole trajectory $X$.
Using the strong Markov property, we get for any function $f$ as above,
\bq
\EE[f(X_t)-\gamma(f)\vert \cF_{\tau^*_M}]
&=&P_{t-\tau^*_M\wedge t}[f](X_{\tau^*_M\wedge t})-\gamma(f)\eq
where $(P_t)_{t\geq 0}$ is the semi-group generated by $L$.
Taking into account that $\sigma(\tau^*_M)$, the $\sigma$-field generated by $\tau^*_M$,
is included into $\cF_{\tau^*_M}$ and that $X_{\tau^*_M}$ is independent from $\tau^*_M$ and distributed according to $\gamma_{[-M,M]}$, we get on the event $\{\tau_M^*\leq t\}$,
\bq
\EE[f(X_t)-\gamma(f)\vert \sigma(\tau^*_M)]&=&\EE[\EE[f(X_t)-\gamma(f)\vert  \cF_{\tau^*_M}]\vert \sigma(\tau^*_M)]\\
&=&\EE[P_{t-\tau^*_M}[f](X_{\tau^*_M})-\gamma(f)\vert \sigma(\tau^*_M)]\\
&=&\EE\lt[
\int P_{t-\tau^*_M}[f](x)\,\gamma_{[-M,M]}(dx)-\gamma(f)\Big\vert \sigma(\tau^*_M)\rt]\\
&\leq & 2
\lVe \nu_{t-\tau^*_M}-\gamma\rVe_{\mathrm{tv}}
\eq
where for any $s\in[0,t]$, $\nu_{t-s}\df\gamma_{[-M,M]}P_{t-s}$ is the law of $X_{t-s}$, when $X$ is started from the initial distribution
$\gamma_{[-M,M]}$.
As a consequence of the Jensen inequality (relatively to the absolute value), it is well-known that the
mapping
\bq
\RR_+\ni s&\mapsto & \lVe \nu_{s}-\gamma\rVe_{\mathrm{tv}}\eq
is non-increasing, so we have proved that
\bq
\EE[f(X_t)-\gamma(f)\vert \sigma(\tau^*_M)]\un_{\{\tau_M^*\leq t\}}&\leq &  2\lVe \gamma_{[-M,M]}-\gamma\rVe_{\mathrm{tv}}\eq
The announced result is a consequence of this bound, by writing
\bq
\EE[f(X_t)-\gamma(f)]&=&\EE[(f(X_t)-\gamma(f))\un_{\{\tau_M^*> t\}}]+\EE[(f(X_t)-\gamma(f))\un_{\{\tau_M^*\leq t\}}]\\
&\leq &2\PP[\tau_M^*> t]+\EE[\EE[f(X_t)-\gamma(f)\vert \sigma(\tau^*_M)]\un_{\{\tau_M^*\leq t\}}]\\
&\leq &2\PP[\tau_M^*> t]+2\lVe \gamma_{[-M,M]}-\gamma\rVe_{\mathrm{tv}}\eq
 and of \eqref{tv2}, by taking the supremum over all  measurable functions $f$ taking values in $[-1,1]$.
\wwtbp\par
The last term of the previous bound is immediate to evaluate:
\begin{lem}
For all $M>0$, we have
\bq
\lVe \gamma_{[-M,M]}-\gamma\rVe_{\mathrm{tv}}&\leq &  \frac{\sqrt{2}}{\sqrt{\pi}M}\exp(-M^2/2)\eq
\end{lem}
\proof
One sees that
\bq
\fo x\in\RR,\qquad \frac{d\gamma_{[-M,M]}}{d\gamma}(x)&=&\frac1{\gamma([-M,M])}\un_{[-M,M]}(x)\eq
so coming back to \eqref{tv}, it appears that
\bq
\lVe \gamma_{[-M,M]}-\gamma\rVe_{\mathrm{tv}}&=&\int_{-M}^M\frac1{\gamma([-M,M])}-1\,d\gamma\\
&=&1-\gamma([-M,M])\\
&=&2\gamma((M,+\iy))\\
&\leq &  \frac{\sqrt{2}}{\sqrt{\pi}M}\exp(-M^2/2)\eq
\wwtbp
\par
In view of Lemma \ref{tvtv}, it remains to study the queues of the distribution of $\tau^*_M$.
The first idea is to use a probabilistic approach via natural comparisons of $Y^*$ with 
simpler processes. This is presented in the appendix, where the weakness of this method is also explained.
Indeed the efficient approach is via spectral considerations in the direction suggested by Remark \ref{Hilbert}.\par\sm
In the above computations, only $Y^*$ was needed, so $X$ (starting from 0) was in fact intertwined with $Y^*$.
It is convenient to adopt the corresponding notations.
Let $L^\dagger$ be the generator of $Y^*$: it acts on functions $f\in \cCc^\iy((0,+\iy))$ 
via
\bq
\fo y\in\RR_+,\qquad L^\dagger[f](y)&\df& f''(y)+V'(y)f'(y)\eq
where
\bq
\fo y\in\RR_+,\qquad V(y)&\df& \frac{y^2}{2}+2\ln(\gamma([0,y]))\eq
So $L^\dagger$ factorizes under the form $\exp(-V)\pa \exp(V)\pa$, making it apparent  that $L^\dagger$
is symmetric in $\LL^2(\nu)$, where $\nu$ is the $\sigma$-finite measure on $\RR$ whose density is $\exp(V)$.
Thus $L^\dagger$ can be extended into its Freidrich extension in $\LL^2(\nu)$.
We will denote $(P^\dagger_t)_{t\geq 0}$ the associated semi-group. 
At least on functions of $\LL^2(\nu)$ which are non-negative, it coincides with its probabilistic
representation  given on measurable and non-negative functions $f$
by
\bq
\fo y\in\RR_+,\qquad P^\dagger_t[f](y)&=&\EE_y[f(Y_t^*)]\eq
where the $y$ in index of the expectation indicates that $Y^*$ starts from $y$.
\\
Besides, we designate by $\Lambda^\dagger$ the Markov kernel from $\RR_+$ to $\RR$ inherited from $\Lambda$: 
\bq
\fo y\in\RR_+,\,\fo A\in\cB(\RR_+),\qquad
\Lambda^\dagger(y,A)&\df&\lt\{\begin{array}{ll}
\delta_{0}(A)&\hbox{, if $y=0$}\\
\noalign{\vskip 2mm}
\frac{\gamma([-y,y]\cap A)}{\gamma([-y,y])}
&\hbox{, otherwise}
\end{array}\rt.\eq
From the previous considerations, we deduce the intertwining relation
\bqn{interdagger}
L^\dagger\Lambda^\dagger&=&\Lambda^\dagger L\eqn
This weak conjugacy relation suggests that the spectral decomposition of $L^\dagger$ should be related to that of $L$.
So let us recall the latter.
Consider $(H_n)_{n\in \ZZ_+}$ the Hermite polynomials defined by
\bq
\fo n\in\ZZ_+,\,\fo x\in\RR,\qquad H_n(x)&\df& (-1)^n\exp(x^2/2)\pa^n\exp(-x^2/2)\eq
They form a orthogonal basis of $\LL^2(\gamma)$ and diagonalize $L$:
\bq
\fo n\in\ZZ_+,\qquad L[H_n]=-nH_n\eq
Note that $H_n$ is even (respectively odd) if $n$ is even (resp.\ odd).
It follows that $\Lambda^\dagger[H_n]=0$ if $n$ is even.
Since $H_0=\un$, we get that $\Lambda^\dagger[H_0]=\un$ and this function does not belong
to $\LL^2(\nu)$ because $\nu$ has an infinite mass.
For the remaining Hermite polynomials, we have:
\begin{lem}\label{H2ndagger}
For all $n\in\NN$, denote 
$H_{2n}^\dagger\df \Lambda^\dagger[H_{2n}]$. This function belongs to $\LL^2(\nu)$,
satisfies $L^\dagger H_{2n}^\dagger =-2nH_{2n}^\dagger $
and is given by
\bq
\fo y>0,\qquad H_{2n}^\dagger(y)&=&-\frac1{\sqrt{2\pi}\gamma([0,y])}H_{2n-1}(y)\exp(-y^2/2)
\eq
\end{lem}
\proof
Indeed, we compute that for any $n\in\NN$ and $y>0$,
\bq
\Lambda^\dagger[H_{2n}](y)&=&\frac1{\sqrt{2\pi}\gamma([-y,y])}
\int_{-y}^yH_{2n}(x)\exp(-x^2/2)\,dx\\
&=&\frac1{\sqrt{2\pi}\gamma([-y,y])}
\int_{-y}^y\pa^{2n}\exp(-x^2/2)\,dx\\
&=&\frac1{\sqrt{2\pi}\gamma([-y,y])}\lt[\pa^{2n-1}\exp(-x^2/2)\rt]_{-y}^{y}\\
&=&\frac2{\sqrt{2\pi}\gamma([-y,y])}\pa^{2n-1}\exp(-y^2/2)\\
&=&-\frac1{\sqrt{2\pi}\gamma([0,y])}H_{2n-1}(y)\exp(-y^2/2)\eq
Thus recalling that for $y>0$, $\nu(y)= (\gamma([0,y]))^2\exp(y^2/2)$,
we get that
\bq
\nu[(H^\dagger_{2n})^2]&=&\frac{1}{2\pi}\int_0^{+\iy} H^2_{2n-1}(y) \exp(-y^2/2)\, dy\\
&=&\frac{1}{\sqrt{2\pi}}\gamma[H^2_{2n-1}]\\
&=&(2n-1)!\eq
(taking into account that for any $n\in\ZZ_+$, $\gamma[H_n^2]=\sqrt{2\pi}n!$).
In particular, $H^\dagger_{2n}$ belongs to $\LL^2(\nu)$ for $n\in\NN$.\par
The fact that $H^\dagger_{2n}$ is an eigenfunction associated to the eigenvalue $-2n$ 
is a consequence of \eqref{interdagger} applied to $H_{2n}$.\wwtbp
\par
Let $\eta$ be the positive measure on $\RR_+$ whose density is given by
\bq
\fo y>0,\qquad \eta(y)&\df& y\gamma([0,y])\eq
It has an infinite weight, but it should nevertheless be seen as a quasi-invariant measure:
\begin{lem}
For all $t\geq 0$ and all measurable and non-negative function $f\st \RR_+\ri\RR_+$,
we have (in $\RR_+\sqcup\{+\iy\}$),
\bq
\eta[P_t^\dagger[f]]&=&\exp(-2t)\eta[f]\eq
\end{lem}
\proof
Consider $H^\dagger_2$, from Lemma \ref{H2ndagger} we have for all $t\geq 0$,
$P^\dagger_t[H^\dagger_2]=\exp(-2t)H^\dagger_2$.
So for any $f\in\LL^2(\nu)$,
\bq
\nu[H^\dagger_2 P_t[f]]&=&\nu[P_t[H^\dagger_2] f]\\
&=&\exp(-2t)\nu[H^\dagger_2 f]\eq
This is the identity announced in the lemma, at least for $f\in\LL^2(\nu)$,
as a consequence of the proportionality of the densities $\eta$ and $\nu H_2^\dagger$:
\bq
\fo y>0,\qquad \nu(y)H_2^\dagger(y)&=&-(\gamma([0,y]))^2\exp(-y^2/2)\frac1{\sqrt{2\pi}\gamma([0,y])}H_{1}(y)\exp(-y^2/2)\\
&=&-\frac1{\sqrt{2\pi}}\gamma([0,y])H_1(y)\\
&=&-\frac1{\sqrt{2\pi}}\eta(y)\eq
The extension to all measurable and non-negative functions $f$ comes 
from the representation of $P_t^\dagger$ as a probability kernel
and from a usual application of the monotone class theorem.
\wwtbp
\par
This result readily shows that the queues of $\tau_M^*$ admits the exponential rate 2,
at least for convenient initial distributions of $Y_0^*$:
\begin{lem}
Assume that the law $m_0$ of $X_0$ has a bounded density with respect to $\eta$. Then there
exists $C>0$ depending on $m_0$ such that
\bq
\fo t\geq 0,\,\fo M>0, \qquad \PP_{m_0}[\tau^*_M>t]&\leq & CM^2\exp(-2t)\eq
\end{lem}
\proof
Denote $f$ the Radon-Nikodym derivate of $m_0$ with respect to $\eta$
and let $C>0$ be an upper bound of $f$.
We have
\bq
 \PP_{m_0}[\tau^*_M>t]&\leq & \PP_{m_0}[Y^*_t\in [0,M]]\\
 &=& \eta[f P_t^\dagger[\un_{[0,M]}]]\\
 &\leq & C\eta [P_t^\dagger[\un_{[0,M]}]]\\
 &\leq &C\exp(-2t)\eta([0,M])\\
 &\leq & C\exp(-2t)\int_0^M y\gamma([0,y])\, dy
 \\
  &\leq & C\exp(-2t)\int_0^M y\, dy
 \\
  &\leq & CM^2\exp(-2t)
 \eq\wwtbp
\par
We want to extend the previous bound to the case where $m_0$ is the Dirac mass at 0.
To do so, first remark that we can restrict ourselves to $M>1$, because $\tau_M^*$ is increasing in $M$.
Next fix $\sigma>0$  small enough such that
$\PP[\tau_1^*\leq \sigma]\leq 1/2$ and denote $\xi$ the sub-probability  which is the image by
$Y^*_\sigma$ of the restriction of $\PP_0$ on $\{\tau_1^*>\sigma\}$.
Its interest is:
\begin{lem}
We have for all $t\geq 0$ and for all $M>1$,
\bq
\PP_{0}[\tau^*_M>\sigma +t]&\leq & 2\PP_\xi[\tau_M^*>t]\eq
\end{lem}
\proof
This is a consequence of the strong Markov property applied to
the stopping time $\sigma\wedge \tau^*_1$:
\bq
\PP_{0}[\tau^*_M>\sigma +t]&=& \EE_0[f(\sigma\wedge \tau_1^*,Y^*_{\sigma\wedge \tau_1^*})]\eq
where \bq
\fo s\in[0,\sigma],\,\fo y\geq 0,\qquad
f(s,y)&\df& \PP_{y}[\tau^*_M>t+\sigma-s]\eq
Note that the quantity $f(s,y)$ is non-decreasing in $s$ and non-increasing in $y$. We deduce
that
\bq
\EE_0[\un_{\{\tau_1^*<\sigma\}}f(\sigma\wedge \tau_1^*,Y^*_{\sigma\wedge \tau_1^*})]&= &
\EE_0[\un_{\{\tau_1^*<\sigma\}}f(\tau_1^*,1)]\\
&\leq & f(\sigma,1)\PP_0[\tau_1^*<\sigma]\eq
Since $\PP_0[\tau_1^*<\sigma]\leq 1/2$ and $f(\sigma,1)\leq f(\sigma,y)$ for all $y\in[0,1]$,
we get
\bq
f(\sigma,1)\PP_0[\tau_1^*<\sigma]&\leq & \EE_0[\un_{\{\tau_1^*\geq \sigma\}}f(\sigma, Y^*_\sigma)]
\\
&\leq & 
\EE_0[\un_{\{\tau_1^*\geq \sigma\}}f(\sigma\wedge \tau_1^*,Y^*_{\sigma\wedge \tau_1^*})]\eq
It follows that
\bq
\PP_{0}[\tau^*_M>\sigma +t]&\leq & 2 \EE_0[\un_{\{\tau_1^*\geq \sigma\}}f(\sigma, Y^*_\sigma)]\\
&=&2 \PP_\xi[\tau_M^*>t]\eq\wwtbp
\par
Thus to prove that there exists a constant $C>0$ such that
\bqn{sigma1}
\fo t\geq \sigma,\,\fo M>1,\qquad  \PP_0[\tau_M^*>t]&\leq & CM^2\exp(-2t)\eqn
it remains to show that $\xi$ admits a density with respect to $\eta$ which is bounded above.
This is not a priori obvious, because $\eta(y)$ is of order $y^2$ for small $y>0$.
But it is true, essentially due to the behavior of the function $g$ defined in \eqref{g} near $0_+$.
\begin{lem}
There exists a constant $C>0$ such that
\bq
\fo y>0 ,\qquad \frac{d\xi}{d\eta}(y)\leq  C\eq
\end{lem}
\proof
Consider the process $Y\df(Y_t)_{t\geq 0}$ starting from 0 and solution of the s.d.e.
\bqn{Bessel3}
\fo t\geq 0,\qquad dY_t&=&\frac2{Y_t}\, dt+\sqrt{2}dB_t\eqn
where $B=(B_t)_{t\geq 0}$ is a standard  Brownian motion. 
Up to the change of time $\RR_+\ni t\mapsto t/2$, $Y$ is a Bessel process of dimension 3.
It follows (see the Section 1 of Chapter 11 of Revuz and Yor \cite{MR1725357}),
that there exists a constant $K>0$ (depending on $\sigma$) such that
the density $\chi$ of $Y_\sigma$ has the form $Ky^2\exp(-y^2/(4\sigma))$.
In particular, we can find another constant $K'>0$ such that 
\bqn{Kprime}
\fo y>0,\qquad \frac{\chi}{\eta}(y)&\leq & K'\eqn
To compare with the law of $Y^*_\sigma$, we use the Girsanov's formula.
More precisely,
define the function $\varphi$ on $\RR_+$ by
\bq
\fo y>0,\qquad \varphi(y)&\df& \frac1{{2}}\int_0^yu+2\frac{\gamma(u)}{\gamma([0,u])}-\frac2{u}\, du\eq
Elementary computations show that these integrals are well-defined, because the integrand is equivalent to $u/3$ for $u>0$ small.
It also appears that $\lve \varphi\rve$, $\lve \varphi'\rve$ and $\lve \varphi''\rve$, as well $u\mapsto \lve \varphi'(u)/u\rve$ are bounded on $(0,\sigma]$.
Since the s.d.e.\ satisfied by $Y^*$ can be written
\bq
\fo t\geq 0,\qquad 
dY^*_t&=& \lt(\frac2{Y_t^*}+\varphi'(Y_t^*)\rt)\,dt +\sqrt{2}dB_t\eq
Girsanov's formula (e.g.\ Chapter 8  of Revuz and Yor \cite{MR1725357}) gives us
\bq
\fo y>0,\qquad
\frac{\xi(y)}{\chi(y)}&=&\EE_0\lt[\exp\lt(\sqrt{2}\int_0^\sigma \varphi'(Y_s)\, dB_s-\int_0^\sigma(\varphi'(Y_s))^2\, ds\rt)\un_{\tau_1(Y)>\sigma}
\Big\vert Y_\sigma=y\rt]\eq
where $Y $ is the solution of \eqref{Bessel3} starting from 0 and $\tau_1(Y)$ is its reaching time of 1.
To evaluate the latter conditional expectation, we write that
\bq
\sqrt{2}\int_0^\sigma \varphi'(Y_s)\, dB_s&=&\varphi(Y_\sigma)
-\int_0^\sigma \frac{2}{Y_s}\varphi'(Y_s)+\varphi''(Y_s)\,ds\eq
which enables to see that
\bq
\lefteqn{\EE_0\lt[\exp\lt(\sqrt{2}\int_0^\sigma \varphi'(Y_s)\, dB_s-\int_0^\sigma(\varphi'(Y_s))^2\, ds\rt)\un_{\tau_1(Y)>\sigma}
\Big\vert Y_\sigma=y\rt]}\\
&=&
\exp(\varphi(y))\EE_0\lt[\exp\lt(-\int_0^\sigma\psi(Y_s)\, ds\rt)\un_{\tau_1(Y)>\sigma}
\Big\vert Y_\sigma=y\rt]
\eq
where
\bq
\fo y>0,\qquad
\psi(y)&\df&  \frac{2}{y}\varphi'(y)+\varphi''(y)+(\varphi'(y))^2\eq
From our previous observations, $\lve \psi(y)\rve$ and  $\lve \varphi(y)\rve$ are bounded for $y\in(0,\sigma]$.
It follows that the function $\xi/\chi$ is bounded on $(0,\sigma]$. This also true on $(\sigma,+\iy)$,
since $\xi$ vanishes there.
In conjunction with \eqref{Kprime}, it ends the proof of the lemma.
\wwtbp\par
Note that there is no difficulty in transforming \eqref{sigma1} into
\bq
\fo t\geq0,\,\fo M>0,\qquad  \PP_0[\tau_M^*>t]&\leq & C(1\vee M)^2\exp(-2t)\eq
up to a change of the constant $C>0$.
Thus putting together all the previous results, we have  proven that there exists a constant $C>0$ such that
for all $t\geq 0$ and all $M>0$, we have
\bq
\lVe m_t-\gamma\rVe_{\mathrm{tv}}&\leq &C(1\vee M)^2\exp(-2t)+
 \frac{\sqrt{2}}{\sqrt{\pi}M}\exp(-M^2/2)\eq
 One could try to minimize the r.h.s.\ in $M>0$ for fixed $t\geq 0$, but it is sufficient to take $M=\sqrt{2t}$
 to see that $\lVe m_t-\gamma\rVe_{\mathrm{tv}}$ converges exponentially fast to zero and that
 \bq
 \limsup_{t\ri+\iy}\frac1{t}\ln(\lVe m_t-\gamma\rVe_{\mathrm{tv}})&\leq &-2\eq
Lemma \ref{taux2} shows that we have recovered the optimal rate, so that the approach via strong times 
is quite sharp.
\par
\begin{rem}
Denote by $\cH$ the Hilbert space generated by  the $H_{2n}$ with $n\in\NN$.
The operator $\Lambda^\dagger$  is compact from $\cH$ to $\LL^2(\nu)$
and one to one.
Indeed, this is an immediate consequence of
\bq
\fo n,m\in\NN,\qquad \nu[H^\dagger_{2n}H^\dagger_{2m}]&=&\frac1{2\sqrt{2\pi}n}\gamma[H_{2n}H_{2m}]\eq
which is shown as in the proof of Lemma \ref{H2ndagger}.
This leads to introduce $\cG\df  \Lambda^\dagger(\cH)$  and to check that $\cG$ is dense in $\LL^2(\nu)$.
It is sufficient to see that any 
smooth mapping $F\st\RR_+\ri \RR$ with compact support belongs to $\Lambda^\dagger(\cH)$, i.e.\ that we can find
a measurable function $f\st (0,+\iy)\ri \RR$ with $\int_0^{+\iy} f\,d\gamma=0$, $\int_0^{+\iy} f^2\,d\gamma<+\iy$
and
\bq
\fo y>0,\qquad \frac{\int_0^y f\,d\gamma}{\gamma([0,y])}&=&F(y)\eq
(we will then have $F=\Lambda^\dagger[\wi f]$ where $\wi f$ is the symmetrization of $f$, which belongs to $\cH$).
So just take
\bq
\fo x>0,\qquad f(x)&\df& \pa( F(x)\gamma([0,x]))\eq
\par\sm
It follows that $(H^\dagger_{2n})_{n\in\NN}$ is an orthogonal Hilbertian basis of $\LL^2(\nu)$ consisting of eigenvectors
of $L^\dagger$. Thus the spectrum of $L^\dagger$ is $-2\NN$. By self-adjointness, we deduce that
\bq
\fo t\geq 0,\,\fo f\in\LL^2(\nu),\qquad
 \lVe P_t[ f]\rVe_{\LL^2(\nu)}&\leq & \exp(-2t)\lVe f\rVe_{\LL^2(\nu)}\eq
This could also have been used to recover the  exponential rate 2 in \eqref{sigma1},
nevertheless we find it more instructive   to work with the quasi-stationary measure $\eta$.
\par\sm
In the same spirit as Remark \ref{Mao}, 
taking into account Theorem 3.3 of the recent preprint of Cheng and Mao \cite{Cheng_eigentime},
we could also have deduced that $Y^*$ is non-explosive from the fact  that sum of the inverse of the eigenvalues
of $-L^\dagger$ in $\LL^2(\nu)$
is infinite, namely, according to the previous considerations, from $\sum_{n\in\NN}1/(2n)=\iy$.
For more informations on the eigentime identity, which states that
 certain reversible Markov processes are explosive if and only if the sum of the inverse of its eigenvalues
is finite, we refer to the paper of Mao
\cite{MR2202590}.
\end{rem}

\appendix
\section{A probabilistic estimate on queues of $\tau^*_M$}
In the previous section we have seen that is important to upper bound
quantities like $\PP[\tau_M^*>t]$ and we obtained nice estimates via
spectral considerations. We were lucky because the spectral decomposition
of $L$ is explicit in the Ornstein-Uhlenbeck example.
In general a probabilistic approach is more flexible, even if in the example at hand
we did not succeed in recovering the optimal rate using this method. Let us nevertheless present 
this approach. At the end we will  see another interplay between probability and spectral theories.\par\me
The basic idea is to compare $Y^*$ with the simpler process $Y\df(Y_t)_{t\geq 0}$ starting from 0
and solution of the s.d.e.
\bqn{sdeY}
\fo t\geq 0,\qquad dY_t&=&Y_t\,dt+\sqrt{2} dB_t\eqn
where $B\df(B_t)_{t\geq 0}$ is a Brownian motion.
We then define for all $M>0$,
\bq
\tau_M&\df& \inf\{t\geq 0\st \lve Y_t\rve =M\}\eq
\begin{lem}
The law of $\tau_M^*$ is stochastically dominated by that of $\tau_M$.
\end{lem}
\proof
Recall  the following behaviors of the mapping $g$ defined in \eqref{g}: as $y$ goes to $0_+$, $g(y)\sim 2/y$
and as $y$ goes to $+\iy$, $g(y)\ll 1/y$.
So we can define
\bq
a&\df& \inf\{ y>0\st g(y)=1/y\}\eq
We first compare $Y^*$ and $Y$ up to the time $\tau^*_a$.
Let $\wi Y$ be an independent  copy of $Y$: it starts from 0 and is solution of the s.d.e.
\bq
\fo t\geq 0,\qquad d\wi Y_t&=&\wi Y_t\,dt+\sqrt{2} d\wi B_t\eq
where $\wi B\df(\wi B_t)_{t\geq 0}$ is a Brownian motion independent from $B$.
Consider the process $\wit Y\df(\wit Y_t)_{t\geq 0}$
given by
\bq
\fo t\geq 0,\qquad \wit Y_t&\df& \sqrt{Y_t^2+\wi Y_t^2}\eq
Simple Itô's computations lead to the fact that $\wit Y$ is the solution starting from 0 of the s.d.e.
\bq
\fo t\geq 0,\qquad d\wi Y_t&=&\lt(\wi Y_t+\frac1{\wi Y_t}\rt)\,dt+\sqrt{2} dW_t\eq
where $W\df(W_t)_{t\geq 0}$ is the Brownian motion defined by
\bq
\fo t\geq 0,\qquad W_t&\df& \int_0^t\frac1{\sqrt{Y_s^2+\wi Y_s^2}}\,(Y_s\, dB_s+\wi Y_s\, d\wi B_s)\eq
Comparing with \eqref{Ystar}, where we replace $B$ with $W$, it appears that 
$\lve Y_t\rve\leq \wit Y_t\leq  Y_t^*$, at least for $t\leq \tau^*_a$.
In particular, for any $M\in(0,a]$, 
the law of $\tau_M$ is stochastically dominated by that of $\tau_M^*$.
Using the strong Markov property at $\tau_M^*$, to prove the same domination
for $M>a$, it is sufficient to deal with the following situation.
Assume that $Y^*$ is the solution of \eqref{Ystar}  starting from $a$ and that $Y$ is solution of \eqref{sdeY} with 
an initial distribution supported by $[0,a]$. Let $B$ be the same in 
\eqref{Ystar}  and in \eqref{sdeY}, then a.s., for all $t\geq 0$, $\lve Y_t\rve\leq Y^*_t$.
\\
Indeed, using Tanaka's formula (see for instance Chapter 6 of the book of Revuz and Yor \cite{MR1725357}),
we have
\bq
\fo t\geq 0,\qquad d\lve Y_t\rve &=&\lve Y_t\rve\,dt+\sqrt{2} dB_t+dl_t\eq
where $(l_t)_{t\geq 0}$ is the local time at 0 of $Y$.
Consider 
\bq
\sigma&\df& \inf\{t\geq 0\st \lve Y_t\rve >Y_t^*\}\eq
If $\sigma <+\iy$, then we have $Y_\sigma=Y^*_\sigma$. 
Recall from Section \ref{dotdp} that necessarily $Y^*_\sigma>0$ and since $l_t$ is only increasing when $Y_t=0$,
there exists a random interval of the form $[\sigma,\sigma+\epsilon)$ on which this local time remains constant.
But we have
\bq
\fo t\geq 0,\qquad
d(Y_t^*-\lve Y_t\rve)&=&(Y_t^*-\lve Y_t\rve)\,dt+ g(Y^*_t)\, dt- dl_t\eq
which, via the parameter variation method,  leads to 
\bq
\fo t\geq 0,\qquad Y_{\sigma +t}^*-\lve Y_{\sigma +t}\rve&=&e^t\int_0^t e^{-s}\, ( g(Y^*_s)\, ds- dl_s)\eq
If $t\in[0,\epsilon)$, the r.h.s.\ is non-negative, which in contradiction with the definition of $\sigma$.\wwtbp
\par
In particular, we get that 
\bqn{weaklink}
\fo M> 0,\,\fo t\geq 0,\qquad \PP[\tau^*_M>t]&\leq & \PP[\tau_M>t]\eqn
The advantage is that the r.h.s. is simpler to evaluate:
\begin{lem}\label{lun}
We have for any $M> 0$ and any $t\geq 0$,
\bq
 \PP[\tau_M>t]&\leq & \sqrt{\frac{2}{(1-e^{-2t})\pi}}e^{-t}M\eq
 \end{lem}
 \proof
 Using once again the parameter variation method, we get that
 \bq
 \fo t\geq 0,\qquad
 Y_t&=& \sqrt{2}\int_0^t\exp(t-s)\, dB_s\eq
 in particular, $Y_t$ is a centered Gaussian random variable of variance $e^{2t}-1$.
 Besides, by definition, we have
 \bq
 \PP[\tau_M>t]&=&\PP[\fo s\in[0,t],\, \lve Y_s\rve\leq  M]\\
 &\leq & \PP[\lve Y_t\rve\leq  M]
 \\
 &=&\int_{-M}^M \exp\lt(-\frac{y^2}{2(e^{2t}-1)}\rt)\, \frac{dy}{\sqrt{2\pi (e^{2t}-1)}}\\
 &\leq & \frac{2M}{\sqrt{2\pi (e^{2t}-1)}}
 \eq
 \wwtbp
\par
These computations leads to the bound
\bq
\fo t>0,\,\fo M>0,\qquad\PP[\tau_M^*>t]
&\leq &  \sqrt{\frac{2}{(1-e^{-2t})\pi}}Me^{-t}\eq
which asymptotically for $t>0$ large, has not the optimal exponential rate (1 instead of 2).\par
So 
where is the weak link in the above arguments?
It is the stochastic dominance \eqref{weaklink}, 
because the exponential rate of $\PP[\tau_M>t]$ for large $t>0$ is almost 1 (for large $M>0$),
as it will be shown below.
So the strong repulsion of $Y^*$ in 0 is the reason for the exponential rate 2  for $\tau_M^*$.
The process $Y$ (or $\lve Y\rve$) has more freedom to wander around 0, which is the best place to ``stay" to avoid the
points $-M$ and $M$, and this accounts for their exit rate 1.\par
\sm
Since the generator of $Y$ is $\wi L\df\exp(-\wi V)\pa \exp(\wi V)\pa$,  where $\wi V\st \RR\ni y\mapsto y^2/2$,
it appears that the measure $\wi\nu$ admitting the density $\exp(\wi V)$
with respect to the Lebesgue measure is ``reversible": $\wi L$ can be extended
into its self-adjoint Friedrich extension on $\LL^2(\wi\nu)$.
From the general Markovian theory of absorption (see e.g.\ the book \cite{MR2986807} of Collet, Mart{\'{\i}}nez and San
        Mart{\'{\i}}n), we have
\bq
\lim_{t\ri+\iy}\frac1{t}\ln(\PP[\tau_M>t])&=&-\lambda_0(M)\eq
where
\bq
\lambda_0(M)&\df& \inf_{f\in\cC^\iy([-M,M])\st f(-M)=f(M)=0}\frac{\int_{-M}^M (f')^2\, d\wi\nu}{\int_{-M}^M f^2\, d\wi\nu}
\eq
Lemma \ref{lun} implies that $\lambda_0(M)\geq 1$ for all $M>0$
and this bound is asymptotically optimal as $M$ goes to infinity:
\begin{lem}
We have
\bq
\lim_{M\ri+\iy} \lambda_0(M)&=&1\eq
\end{lem}
\proof
Let $f_M$ be the function defined on $[-M,M]$ by
\bq
\fo y\in[-M,M],\qquad f_M(y)&\df& \exp(-y^2/2)-\exp(-M^2/2)\eq
Elementary computations show that
\bq
\liminf_{M\ri+\iy}\lambda_0(M)&\geq & \lim_{M\ri+\iy}\frac{\int_{-M}^M (f_M')^2\, d\wi\nu}{\int_{-M}^M f_M^2\, d\wi\nu}\\
&=&1\eq\wwtbp
\par
The functions $f_M$, for $M>0$, were suggested by the spectral decomposition of $\wi L$ on $\LL^2(\gamma)$,
which can be obtained by a method somewhat dual to the one presented in the previous section.
Consider, on the appropriate domain of $\LL^2(\gamma)$, the linear mapping
$K\st f\mapsto \exp(-V)\pa f\in\LL^2(\wi \nu)$.
Since $\wi L=\exp(-\wi V)\pa \exp(\wi V)\pa$ and $L=\exp(\wi V)\pa \exp(-\wi V)\pa$, we get at once
the intertwining relation $\wi LK=KL$ (with a non-Markovian link $K$, but its inverse is a positive kernel
 quite close to $\Lambda^\dagger$).
So a priori  the $\wi H_n\df K[H_n]$, for $n\in\NN$, are good candidates to be the eigenvectors of $\wi L$, associated respectively to the eigenvalues $-n$.
Indeed, we compute that 
\bq
\fo n\in\NN,\,\fo y\in\RR,\qquad \wi H_n(y)&=&n\exp(-y^2/2)H_{n-1}(y)\eq
so that $(\wi H_n)_{n\in\NN}$ is an orthogonal Hilbertian basis of $\LL^2(\wi\nu)$, so the spectrum of $\wi L$ is $-\NN$.
The measure $\wi\eta\df \wi H_1\wi\nu=\exp(-\wi V)$ is quasi-stationary for $\wi L$ and its adaptation to the Dirichlet boundary conditions on $[-M,M]$
furnishes $f_M$, for $M>0$.
One can also deduce the spectral decomposition of the generator of $\lve Y\rve$: restrict everything to $\RR_+$,
but just keep the $\wi H_n$ with $n$ odd. In particular its spectrum is $\{-1,-3, -5, ...\}$.

\bigskip\bigskip
\par\hskip5mm\textbf{\large Aknowledgments:}\par\sm\noindent 
I'm grateful to Li-Juan Cheng and Yong-Hua Mao
for having sent me their unpublished papers \cite{Cheng_eigentime} and 
\cite{Cheng_passage}.

% \bibliography{MR3s,arXiv,nonpubliees,MR-Miclo}
% \bibliographystyle{plain}

\def\cprime{$'$}

\vskip2cm
\hskip70mm\box5

\end{document}